\documentclass[11pt]{article}

\usepackage{latexsym,amssymb,amsmath,graphicx, amsthm, times}
\usepackage{times}

\newcommand {\ME}{\mathbb{E}^{x}}

\renewcommand{\S}{\mathcal{S}}

\numberwithin{equation}{section}
\newtheorem{proposition}{Proposition}[section]

\newtheorem{remark}{Remark}[section]
\newtheorem{lemma}[proposition]{Lemma}

\newtheorem{assump}{Assumption}[section]

\newcommand {\R}{\mathbb{R}}

\newcommand {\F}{\mathcal{F}}
\newcommand {\A}{\mathcal{A}}

\newcommand {\p}{\mathbb{P}}

\newcommand {\E}{\mathbb{E}}
\title{The Effects of Implementation Delay on Decision-Making Under
Uncertainty}
\author{Erhan Bayraktar \thanks{E. Bayraktar was supported in part by the National Science Foundation, under grant
DMS-0604491.}
\and Masahiko Egami
\thanks{E. Bayraktar and M. Egami are in the Department of Mathematics, University
of Michigan, Ann Arbor, MI 48109, USA, email: \{erhan,
egami\}@umich.edu. } }
\date{}
\begin{document}

\maketitle
\begin{abstract}\noindent
In this paper, we accomplish two objectives:  First, we provide a
new mathematical characterization of the value function for
impulse control problems with implementation delay and present a
direct solution method that differs from its counterparts that use
quasi-variational inequalities. Our method is direct, in the sense
that we do not have to guess the form of the solution and we do
not have to prove that the conjectured solution satisfies
conditions of a verification lemma. Second, by employing this
direct solution method, we solve two examples that involve
decision delays: an exchange rate intervention problem and a
problem of labor force optimization.
\end{abstract}
\section{Introduction}

Implementation delays occur naturally in decision-making problems.
Many corporations face regulatory delays, which need to be taken
into account when the corporations make decisions under
uncertainty. A decision made will be carried out only after
certain amount of time elapses, for example,  due to regulatory
reasons.  The decision involves optimally exercising a \emph{real
option} or optimally manipulating (with some associated cost) a
state variable, which is the source of uncertainty. Several
problems that fit into this framework can be found in the
literature: The work of Bar-Ilan and Strange
\cite{bar-ilan-strange1} constitutes the first study considering
how delays affect rational investment behavior. Keppo and Peura
\cite{keppo-peura}
 consider the decision making
problem a bank has to solve when it is faced with a minimum
capital requirement, a random income,  and delayed (and costly)
recapitalization. The bank's problem is to determine when to raise
capital from its shareholders and the amount to be raised, given
that this transaction requires a heavy preparatory work, which
causes delay. Bar-Ilan and Strange \cite{bar-ilan} consider
(irreversible) sequential (2 stage) investment decision problems
given two sources of delay: one due to market analysis in the
first stage and the other due to construction of a production
facility in the second stage.  In each stage the firm's problem is
to decide whether to continue entering into the market (of that
product) or to abandon it. See also Subramanian and Jarrow
\cite{sub-jarrow} who consider the problem of a trader (who is not
a \emph{price taker}) who wants to liquidate her position and
encounters execution delays in an illiquid market.  Alvarez and
Keppo \cite{alvarez-keppo}  study the impact of delivery lags on
irreversible investment demand under revenue uncertainty.
{\O}ksendal et. al. \cite{oksendal0sulem2001},
\cite{elsanosi-oks-sulem} consider the classical stochastic
control of stochastic delays systems.

The problem of finding an optimal decision (in the presence of
delays) can be characterized as a stochastic impulse control
problem or an optimal stopping problem. In the papers cited above
the impulse control problem or the optimal stopping problem were
solved by using a system of quasi-variational inequalities. (See
e.g. Bensoussan and Lions \cite{ben-lions} and {\O}ksendal and
Sulem \cite{oksendal-book-2} for the relationship between control
problems and quasi-variational inequalities.)  In a  different
approach,  {\O}ksendal and Sulem \cite{oksendal-delay-impulse}
solve a version of delay problems, in which the controller decides
on the magnitude of control at the time of decision-making before
any delay (the decision is implemented after some delay). They
convert the optimal impulse control problem with delayed reaction
into a no-delay optimal stopping/impulse control problem. Note
that choosing the control in this way introduces strong path
dependence of the controlled process.

Here, we solve the impulse control problems with delays
\emph{directly} and the magnitude of the impulses are chosen at
the time of action, not at the time of decision-making, by
providing a new characterization of the \emph{value function}. The
controlled process is a non-Markov process in this case, too,
since depending on when a point in the state space is reached, it
has different roles. But the controlled process in this case
regenerates after a decision is implemented, and the value of the
state process during the delay time depends on the past only
through the value of the state process at the time of
decision-making. We will only consider the threshold and band
policies in this paper, since we expect that the non-Markovian
structure will make finding the optimal solution much more
difficult if we allow more general strategies. For example,
because of the lack of Strong Markov property, we were unable to
prove the concavity properties of the value function when the
admissible strategies were a superset of band or threshold
strategies.

Our results rely on the works of Dynkin \cite{dynkin2},
\cite{dynkin} (see e.g. Theorem 16.4) and Dayanik and Karatzas
\cite{DK2003}, who give a general characterization of optimal
stopping times of one dimensional diffusions, and on the work of
Dayanik and Egami \cite{DE2005}, who characterize the value
function of stochastic impulse control problems. Our method is
direct, in the sense that we do not have to guess the form of the
solution and we do not have to prove that the conjectured solution
satisfies conditions of a verification lemma as all the methods in
the above literature do. Other works similar in vein to ours that
 provide different characterizations of the value function of
impulse/singular control problems for one dimensional diffusions
rather than solving variational inequalities are Alvarez
\cite{alvarez2}, \cite{alvarez1}; Alvarez and Virtanen
\cite{alvarez3}; and Weerasighe \cite{wee}.

We give a geometric characterization of the value function,
specifically, we find very general conditions on the reward
function and the coefficients of the underlying diffusion under
which the value function can be linearized (in the continuation
region) after a suitable transformation. Then the problem of
determining the value function is equivalent to determining the
slope (if admissible strategies are threshold strategies), the
slope and the intercept (if admissible strategies are band
strategies) from first order conditions. To show the efficacy of
our methodology we apply it to an optimization problem of a
central bank that needs to carry out exchange rate intervention
(this is the Krugman model of interest rates considered, among
others, in Mundaca and {\O}ksendal \cite{mundaca-oksendal}) when
there is delay in the implementation of its decisions.
 Also, using our methodology we will find optimal hiring
and firing decisions of a firm that faces stochastic demand and
has to conform to regulatory delays. Other works that deal with
 labor optimization problems are Bentolila and Bertola \cite{bb},
 and Shepp and Shiryaev \cite{shepp-shiryaev} who model firing and
 hiring decisions as singular controls. It is also worth pointing
 out that an impulse control study when the underlying process is
 a superposition of a Brownian motion and a compound Poisson
 process (when the jumps are of phase type) is given by
 \cite{bar-ilan-perry-stadje} with management of foreign exchange
 reserves and labor optimization in mind.

The rest of the paper is organized as follows: In Section 2, we
give a characterization of general threshold strategies with
implementation delays and provide an easily implemented algorithm
to find the value function and the optimal control. To illustrate
our methodology, we will solve a delayed version of an example
from Mundaca and {\O}ksendal \cite{mundaca-oksendal} (also see
{\O}ksendal \cite{O1999}).  A similar problem to the one we
consider was solved in \cite{oksendal-delay-impulse} in which the
controller decides on the magnitude of control at the time of
decision-making before any delay. In Section 3, we work with a
band policy. In this section we work on the specific example of
optimal hiring and firing decisions rather than providing a
general characterization for the value function. We again provide
an easily implemented algorithm to find the optimal control.
Finally, we conclude in Section~\ref{sec:conclusion}.

\section{Optimal Threshold Strategies}
 Let $(\Omega,
\F, \p)$ be a complete probability space with a standard Brownian
motion $W=\{W_t; t\geq 0\}$ and consider the diffusion process $X^0$
with state pace $\mathcal{I}=(c, d)\subseteq \mathbb{R}$ and
dynamics
\begin{equation}\label{eq:process}
dX^0_t=\mu(X^0_t)dt + \sigma(X^0_t)dW_t
\end{equation}
for some Borel functions $\mu :\mathcal{I}\rightarrow \mathbb{R}$
and $\sigma :\mathcal{I}\rightarrow (0, \infty)$. (We assume that
the functions $\mu$ and $\sigma$ are sufficiently regular so that
(\ref{eq:process}) makes sense.) Here we take $c$ and $d$ to be a
natural boundaries. We use ``0" as the superscript to indicate
that $X^0$ is the uncontrolled process. We denote the
infinitesimal generator of $X^0$ by $\A$ and consider the ODE
$(\A-\alpha)v(x)=0$. This equation has two fundamental solutions,
$\psi(\cdot)$ and $\varphi(\cdot)$. We set $\psi(\cdot)$ to be the
increasing and $\varphi(\cdot)$ to be the decreasing solution.
$\psi(c+)=0, \varphi(c+)=\infty$ and $\psi(d-)=\infty,
\varphi(d-)=0$ because both $c$ and $d$ are natural boundaries.
First, we define an increasing function
\begin{equation}\label{eq:F}
F(x)\triangleq\frac{\psi(x)}{\varphi(x)}.
\end{equation}
Next, following \cite{dynkin}, p. 238, we define concavity of a
function with respect $F$ as follows: A real valued function $u$
is called \emph{$F$-concave} on $(c, d)$ if, for every $c<l<r< d$
and $x\in[l, r]$,
\begin{equation*}
u(x)\geq
u(l)\frac{F(r)-F(x)}{F(r)-F(l)}+u(r)\frac{F(x)-F(l)}{F(r)-F(l)}.
\end{equation*}

Suppose that at any time $t\in\mathbb{R_+}$ and any state
$x\in\mathbb{R_+}$, we can intervene and give the system an impulse
$\xi\in \mathbb{R}$. Once the system gets intervened, the point
moves from $x$ to $y\in\mathbb{R_+}$ with associated reward and cost.
An impulse control for the system is a double sequence,
\begin{equation}
\nu =(T_1, T_2,....T_i....; \xi_1,\xi_2,...\xi_i....)
\end{equation}
where $0\leq T_1<T_2<....$ is an increasing sequence of
$\mathbb{F}$-stopping times such that $T_{i+1}-T_i \geq \Delta$, and $\xi_1$, $\xi_2...$ are
$\mathcal{F}_{(T_i+\Delta)-}$ measurable random variables
representing impulses exercised at the corresponding intervention
times $T_i$ with $\xi_i\in Z$ for all $i$ where $Z\subset
\mathbb{R}$ is a given set of admissible impulse values.  The
controlled process until the first intervention time is described as follows:
\begin{align} \label{eq:control}
\begin{cases}
dX_t = \mu(X_t)dt + \sigma(X_t)dW_t, \quad 0 \leq t<T_1+\Delta \\
X_{T_{1}+\Delta}=\Gamma(X_{(T_{1}+\Delta)-}, \xi_1)
\end{cases}
\end{align}
with some mapping $\Gamma: (c,d)\times\mathbb{R}\rightarrow
\mathbb{R}$.  We consider the following performance measure
associated with $\nu\in\mathcal{V}$ (= a collection of admissible
strategies),
\begin{equation}\label{eq:J}
    J^\nu(x)=\ME\left[\int_0^{\infty}e^{-\alpha s}f(X_s)ds+\sum_{T_i<\infty}e^{-\alpha
    (T_i+\Delta)}K(X_{(T_i+\Delta)-},X_{T_i+\Delta})\right].
\end{equation}
The objective (we shall call it the \emph{delay problem}) is to
find the optimal strategy $\nu^*$ (if it exists) and the value
function:
\begin{equation}\label{problem1}
v(x)\triangleq \sup_{\nu\in \mathcal{V}}J^{\nu}(x)=J^{\nu^*}(x).
\end{equation}

\begin{remark}\label{rem:not-Markov}
The controlled process $X$ is not a Markov process, since
depending on whether a point is reached in the time interval
$[T_i,T_i+\Delta)$ or not,  that point has different roles. (The controlled process might
jump or not at a given point depending on how it reaches to that
point.) However, 1) the process regenerates at
times $\{T_i+\Delta\}_{i \in \mathbb{N}}$, and 2) the value of the
process at time $T \in (T_i,T_{i+\Delta})$, $X_T$, depends on the
information up to $T_i$, $\mathcal{F}_{T_i}$, only through the
value of the process at time $T_i$, $X_{T_i}$. Instead of finding
the optimal strategy for a non-Markov process, we will use the
hints of Markovian features to find the optimal \emph{threshold
strategy} (see Assumption~\ref{ass:first-assump}).
\end{remark}
The following is a standing assumption in Sections 2.1 and 2.2.
\begin{assump}\normalfont \label{ass:first-assump}
We make the following assumptions in this section:
\begin{enumerate}
\item[(a)] We will assume that the set of admissible strategies is limited to \emph{threshold
strategies}. These strategies are determined by specifying two
numbers $a \in (c,d)$ and $b \in (c, d)$ as follows: At the time
the uncontrolled process hits level $b$, the controller decides to
reduce the level of the process from $\xi_{T_i-}=b$ to $a<b$,
through an intervention, and save the continuously incurred cost
(which is high if the process is at a high level).  But the
implementation of this decision is subject to a delay of $\Delta$
units of time. Note that $\xi_{(T_i+\Delta)-}$ might be less than
$a$. In that case the impulse applied increases the value of the
process. Otherwise, if the value of the process is greater than
$a$ at time $(T_i+\Delta)-$ then the intervention reduces the
level of the process to $a$.

\item [(b)] The running cost function $f:(c,d)\rightarrow \mathbb{R}$ is a
continuous functions that satisfies
\begin{equation} \label{eq:f-condition}
\ME\left[\int_0^\infty e^{-\alpha s}|f(X_s)|ds\right]<\infty.
\end{equation}
\item[(c)]For any point $x\in(c,d)$,  we assume
\begin{equation}\label{eq:fixed-cost}
K(x, x)<0.
\end{equation}
We make this assumption to account for the fixed cost of making an
intervention.
\end{enumerate}
\end{assump}

\subsection{Characterization of the Value
Function}\label{sec:characterization}

In this section, we will show that when we apply a suitable
transformation to the value function corresponding to a particular
threshold strategy (that is identified by a pair $(a,b)$), the
transformed value function is linear on $(0, F(b))$. This
characterization will become important in determining the optimal
threshold strategy in the next section.

Let us define
\begin{equation}
 g(x)\triangleq \ME \left[\int_0^\infty e^{-\alpha
s}f(X_s^0)ds\right]
\end{equation}
The following identity, which can be
derived using the Strong Markov Property of $X^0$, will come handy
in a couple of computations below:
\begin{equation}
\ME \left[\int_{0}^{\tau} e^{-\alpha s} f(X^0_s)
ds\right]=g(x)-\ME\left[e^{-\alpha \tau} g(X^0_{\tau})\right],
\end{equation}
for any stopping time $\tau$ under the assumption
(\ref{eq:f-condition}).

\noindent Now, let us simplify $J^\nu$ by splitting the terms in
(\ref{eq:J}).  We can write the first terms (the term with the
integral) as
\begin{equation}\label{eq:integral-of-f}
\begin{split}
&\ME\left[\int_0^{\infty}e^{-\alpha s}f(X_s)ds\right]\\
&\hspace{1.5cm}=\ME\left[\int_0^{T_1+\Delta}e^{-\alpha
s}f(X^0_s)ds+ e^{-\alpha (T_1+\Delta)}\E^{X_{T_1+\Delta}}
\left[\int_{0}^{\infty}e^{-\alpha s}f(X_s)ds\right]\right]\\
&\hspace{1.5cm}=g(x)-\ME[e^{-\alpha(T_1+\Delta)}g(X^0_{T_1+\Delta})]+\ME\left[e^{-\alpha
(T_1+\Delta)}\E^{X_{T_1+\Delta}}\left[ \int_{0}^{\infty}e^{-\alpha
s}f(X_s)ds\right]\right]\\
&\hspace{1.5cm}=g(x)-\ME[e^{-\alpha(T_1+\Delta)}g(X_{(T_1+\Delta)-})]+\ME\left[e^{-\alpha
(T_1+\Delta)}\E^{X_{T_1+\Delta}}\left[ \int_{0}^{\infty}e^{-\alpha
s}f(X_s)ds\right]\right],
\end{split}
\end{equation}
while the second term can be developed as {\small
\begin{align*}
&\ME\left[\sum_{T_i<\infty}e^{-\alpha
    (T_i+\Delta)}K(X_{(T_1+\Delta)-},X_{T_1+\Delta})\right]\\
&= \ME\left[e^{-\alpha
(T_1+\Delta)}K(X_{(T_1+\Delta)-},X_{T_1+\Delta})+e^{-\alpha(T_1+\Delta)}\sum_{i=2}^{\infty}e^{-\alpha((T_i+\Delta)-(T_1+\Delta))}K(X_{(T_i+\Delta)-},X_{T_i+\Delta})\right]\\
&=\ME\left[e^{-\alpha
(T_1+\Delta)}K(X_{(T_1+\Delta)-},X_{T_1+\Delta})+e^{-\alpha(T_1+\Delta)}\ME\left
[\sum_{i=1}^{\infty}e^{-\alpha((T_i+\Delta)\circ\theta(T_1+\Delta))}
K(X_{(T_{i+1}+\Delta)-},X_{T_{i+1}+\Delta})\bigg|\F_{T_1+\Delta}\right]\right]\\
&= \ME\left[e^{-\alpha
(T_1+\Delta)}\left\{K(X_{(T_1+\Delta)-},X_{T_1+\Delta})+
\E^{X_{T_1+\Delta}}\left[\sum_{i=1}^{\infty}e^{-\alpha(T_i+\Delta)}K(X_{(T_i+\Delta)-},X_{T_i+\Delta})\right]\right\}\right]
\end{align*}}
where we used $T_{i+1}+\Delta=
(T_1+\Delta)+(T_i+\Delta)\circ\theta(T_1+\Delta)$  with the shift
operator $\theta(\cdot)$ in the second equality. Here, we relied
on Remark~\ref{rem:not-Markov}. Combining the two terms, we can
write (\ref{eq:J}) as
\begin{align*}
J^\nu(x)=\ME\left[e^{-\alpha
(T_1+\Delta)}\left\{K(X_{(T_i+\Delta)-},X_{T_i+\Delta})-g(X_{(T_1+\Delta)-})+J^\nu(X_{T_1+\Delta})\right\}\right]+g(x).
\end{align*}
We define
\begin{equation}\label{eq:defn-u-J-g}
u\triangleq J^\nu-g.
\end{equation}
By adding and subtracting $g(X_{(T_1+\Delta)})$ to and from the
first term we obtain
\begin{equation}\label{eq:u1}
u(x)=\ME\left[e^{-\alpha (T_1+\Delta)}\bar{K}(X_{(T_1+\Delta)-},
X_{T_1+\Delta})+u(X_{T_1+\Delta})\right]
\end{equation}
in which
\begin{equation}\label{eq:u1-under}
\bar{K}(x, y)\triangleq K(x, y)-g (x)+g(y).
\end{equation}
 Since $T_1-=\tau_b$ with $\tau_b=\inf\{t\geq 0: X^0_t\geq b\}$
and the post intervention point by
\begin{equation}
X_{T_1+\Delta}=X_{\tau_b+\Delta}=X_{(\tau_b+\Delta)-}-\xi_1\triangleq
a.
\end{equation}
From Remark~\ref{rem:not-Markov}
\begin{align}\label{eq:u2}
u(x)&=\ME\left[e^{-\alpha
(\tau_b+\Delta)}\left\{\bar{K}(X_{\tau_b+\Delta},
a)+u(a)\right\}\right]\nonumber\\
      &=\ME\left[\ME\left[e^{-\alpha
      (\tau_b+\Delta)}\left\{\bar{K}(X_{\tau_b+\Delta},
a)+u(a)\right\}\bigg| \F_{\tau_b}\right]\right]\nonumber\\
      &=\ME\left[e^{-\alpha\tau_b}\E^{X_{\tau_b}}\left[e^{-\alpha\Delta}\left\{\bar{K}(X_{\Delta},
a)+u(a)\right\}\right]\right].
\end{align}
Evaluating at $x=b$, we obtain
$u(b)=\E^b[e^{-\alpha\Delta}\left\{\bar{K}(X_{\Delta},
a)+u(a)\right\}]$. Therefore,  (\ref{eq:u1}) becomes
\begin{equation*}
u(x)=\ME\left[e^{-\alpha\tau_b}u(X_{\tau_b})\right].
\end{equation*}Hence we have finally
\begin{eqnarray}\label{eq:twoside}
u(x)&=
\begin{cases}
u_0(x)\triangleq \ME\left[e^{-\alpha\tau_b}u(b)\right], &x\in(c, b),\\
\ME\left[e^{-\alpha\Delta}(\bar{K}(X_{\Delta}, a)+u_0(a))\right],
&x\in[b, d),
\end{cases}
\end{eqnarray}
where the second equality is obtained when we plug $T_1=0$ in
(\ref{eq:u1}).

\noindent  Using appropriate boundary conditions one can solve
$(\mathcal{A}-\alpha)u=0$ and obtain
\begin{equation}\label{eq:laplace}
\ME[e^{-\alpha\tau_r}1_{\{\tau_r<\tau_l\}}]=\frac{\psi(l)\varphi(x)-\psi(x)\varphi(l)}
{\psi(l)\varphi(r)-\psi(r)\varphi(l)}, \, \,
\ME[e^{-\alpha\tau_r}1_{\{\tau_l<\tau_r\}}]=\frac{\psi(x)\varphi(r)-\psi(r)\varphi(x)}
{\psi(l)\varphi(r)-\psi(r)\varphi(l)},
\end{equation}
for $x\in[l,r]$ where $\tau_l\triangleq\inf\{t>0; X^0_t=l\}$ and
$\tau_r\triangleq\inf\{t>0; X^0_t=r\}$ (see e.g. Dayanik and
Karatzas \cite{DK2003}). By defining
\begin{equation}\label{eq:W-u-trnasfrom}
W \triangleq (u/\varphi)\circ F^{-1},
\end{equation}
equation (\ref{eq:twoside}) becomes
\begin{align}\label{eq:W0}
W(F(x))&=W(F(c))\frac{F(b)-F(x)}{F(b)-F(c)}+W(F(b))\frac{F(x)-F(c)}{F(b)-F(c)},
\quad x\in(c, b],
\end{align}
We should note that $F(c)\triangleq F(c+)=\psi(c+)/\varphi(c+)=0$
and
\begin{equation}
W(F(c))=l_c\triangleq \limsup_{x\downarrow c}\frac{\bar{K}(x,
a)^+}{\varphi(x)}
\end{equation}
for any $a\in (c,d)$.  For more detailed mathematical meaning of
this value $l_c$, we refer the reader to Dayanik and
Karatzas\cite{DK2003}.  We have now established that $W(F(x))$ is
a \emph{linear function} in the transformed ``continuation
region".

\subsection{An Algorithm to Compute the Value
Function}\label{sec:algorithm-threshold}
Let us denote
\begin{equation}\label{eq:defn-r-forex}
r(x;a)\triangleq \ME[e^{-\alpha\Delta}\bar{K}(X_\Delta, a)]
\end{equation}
and transform this function by
\begin{equation}\label{eq:R-forex}
R(\cdot; a)\triangleq \frac{r(F^{-1}(\cdot),
a)}{\varphi(F^{-1}(\cdot))}.
\end{equation}

\noindent \underline{First stage}: For a given pair $(a,b) \in
(c,d) \times (c,d)$ we can determine (\ref{eq:twoside}) from the
linear characterization (\ref{eq:W0}). On $(0, F(b)]$ we will find
$W(y)=\rho y+ l_c$ (in which the slope is to be determined) from
\begin{equation}\label{eq:cont-fit-FoREX}
\rho F(b)+l_c=R(F(b), a)+e^{-\alpha \Delta}(\rho
F(a)+l_c)\frac{\varphi(a)}{\varphi(b)}.
\end{equation}
$\rho$ can be determined as
\begin{equation}\label{eq:rho-FOREX}
\rho=\frac{R(F(b;a))+l_c(e^{-\alpha \Delta}\frac{\varphi(a)}{\varphi(b)}-1)}{F(b)
-e^{-\alpha \Delta} \frac{\varphi(a)}{\varphi(b)}F(a)}
\end{equation}
Sometimes we will refer to $\rho$ as $b \rightarrow \rho(b)$, when
it becomes necessary to emphasize the dependence on $b$. The
function $u$ can be written as
\begin{equation}\label{eq:u-lb-gb}
u(x)=\begin{cases} u_{0}(x) \triangleq \rho \psi(x)+l_c \varphi(x)
& x \leq b
\\
r(x,a)+e^{-\alpha \Delta}u_0(a) & x>b.
\end{cases}
\end{equation}
Note that $(\A-\alpha)u(x)=0$ for $x<b$. \emph{Henceforth, to
emphasize the dependence on the pair $(a,b)$ we will write
$u^{a,b}(\cdot)$ for the function $u(\cdot)$}.

\noindent \underline{Second stage}: Our purpose in this section is
to determine
\begin{equation}
u^{a}(x)\triangleq \sup_{b \in (c,d)}u^{a,b}(x), \quad x \in
(c,d),
\end{equation}
to determine the constant $b^*$
\begin{equation}\label{eq:max-attained}
u^{a}(x)=u^{a,b^*}(x), \quad x\in (c,d),
\end{equation}
if there exists one.

Let us fix $a$ and treat $\rho$ as a function of $b$ parametrized
by a.
\begin{lemma}
Assume that the function $R(\cdot; a)$ defined in
(\ref{eq:R-forex}) is differentiable and that there exists a
constant $b^* \in (c,d)$ satisfying (\ref{eq:max-attained}). Then
$b^*$ satisfies the equation
\begin{equation}\label{eq:non-linear-ipmlicit-FOREX}
\rho F'(b)=\frac{\partial}{\partial y} R(y;a)
\bigg|_{y=F(b)}F'(b)- e^{-\alpha \Delta}(\rho
F(a)+l_c)\frac{\varphi(a) \varphi'(b)}{\varphi(b)^2}.
\end{equation}
in which $\rho$ is given by (\ref{eq:rho-FOREX}).
\end{lemma}
\begin{proof}
From (\ref{eq:u-lb-gb}) it follows that the maximums of the
functions $b \rightarrow u^{a,b}$ and $b \rightarrow \rho(b)$ are
attained at the same point. Now taking the derivative of
(\ref{eq:cont-fit-FoREX}) and evaluating at $\rho_b=0$ we obtain
(\ref{eq:non-linear-ipmlicit-FOREX}).
\end{proof}
 To find the
optimal $b$ (given $a$) we solve the non-linear and implicit
equation (\ref{eq:non-linear-ipmlicit-FOREX}). Under certain
assumptions on the function $(r/\varphi)\circ F^{-1}$, this
equation has a unique solution as we show below.

\begin{remark}\label{rem:smooth-fit-forex}
On $y \geq F(b)$, the function $W$ is given by
\begin{equation}\label{eq:rem1}
W(y)=e^{-\alpha \Delta}(\rho
F(a)+l_c)\frac{\varphi(a)}{\varphi( F^{-1}(y))}+R(y;a).
\end{equation}
The right derivative of $W$ at $F(b)$ is given by
\begin{equation}\label{eq:rem2}
W'(F(b))=-e^{-\alpha \Delta}(\rho
F(a)+l_c)\frac{\varphi(a)}{\varphi(b)^2}\frac{\varphi'(b)}{F'(b)}+\frac{\partial}{\partial
y}R(y;a)\bigg|_{y=F(b)}.
\end{equation}
Therefore, (\ref{eq:non-linear-ipmlicit-FOREX}) implies that the
 left
and the right derivative of $W$ (recall that $W(y)=\rho y+l_c$ for $y<F(b)$) at $F(b)$ are equal (smooth fit).
\end{remark}

Let us define
\begin{equation} \label{eq:inner-optimization}
  u_a(x)\triangleq \sup_{b \in (c,d)}\ME\left[e^{-\alpha\tau_b}\E^{X_{\tau_b}}\left[e^{-\alpha\Delta}\left\{\bar{K}(X_{\Delta},
a)+u_a(a)\right\}\right]\right].
\end{equation}
The next lemma shows that (\ref{eq:inner-optimization}) is
well-defined. Below we show that under certain assumptions on
$(r/\varphi)\circ F^{-1}$ this function is equal to $u^a$.
\begin{lemma}\label{lem:existence-lemma}
Assume that
\begin{equation}\label{eq:positivity}
\sup_{x\in (c,d)}\ME[\bar{K}(X_\Delta, a)]>0
\end{equation}
for some $a\in (c,d)$. Let us introduce a family of value
functions parameterized by $\gamma \in \mathbb{R}$ as
\begin{equation}\label{eq:param}
V_a^\gamma(x)\triangleq\sup_{\tau\in
\S}\ME\left[e^{-\alpha(\tau+\Delta)}\left\{\bar{K}(X^0_{\tau+\Delta},
a)+\gamma\right\}\right]=\sup_{\tau\in\S}\ME\Big[e^{-\alpha\tau}\E^{X^0_{\tau}}\left[e^{-\alpha\Delta}\left\{\bar{K}(X^0_\Delta,
a)+\gamma\right\}\right] \Big],
\end{equation}
here $\mathcal{S}$ is the set of all stopping times of the
filtration natural filtration of $X^0$. Then there exists a unique
$\gamma^*$ such that $V_a^{\gamma^*}(a)=\gamma^*$.
\end{lemma}
\begin{proof}
 Let us denote
\begin{equation}\label{eq:Wa} W^\gamma_a(F(x)) \triangleq
\frac{V^\gamma_a(x)}{\varphi(x)},
\end{equation}
Consider the function $\gamma \rightarrow V^\gamma_a(a)$. Our aim
is to show that there exists a fixed point to this function. Let
us consider $V^0_a(a)$ first. Because (\ref{eq:positivity}) is
satisfied we have that  $V^0_a(a)>0$.
 As $\gamma$ increases,
$V^\gamma(a)$ increases monotonically, by the right hand side of
(\ref{eq:param}). Now, Lemma~\ref{lem:Lipschitz} implies that for
$\gamma_1>\gamma_2\geq 0$,
\begin{equation}\label{eq:contraction2}
V_a^{\gamma_1}(x)-V_a^{\gamma_2}(x)\leq \gamma_1-\gamma_2
\end{equation}
for any $x\in\mathbb{R_+}$. Note that  $W^\gamma_a(F(a)) \geq
R(F(a), a)+\frac{e^{-\alpha\Delta}\gamma}{\varphi(a)}$ for all
$\gamma$. However, since $V$ has less than linear growth in
$\gamma$ as demonstrated by (\ref{eq:contraction2}) we can see
that
 there is a certain $\gamma^{'}$ large
enough such that $W^\gamma_a(F(a))=R(F(a),
a)+\frac{e^{-\alpha\Delta}\gamma}{\varphi(a)}$ for $\gamma \geq
\gamma^{'}$. This implies however
\begin{align*}
\varphi(a)W^{\gamma^{'}}_a(F(a))&=\varphi(a)R(F(a),
a)+e^{-\alpha\Delta}\gamma'\\ \Leftrightarrow
V_a^{\gamma'}(a)&=r(a, a)+e^{-\alpha\Delta}\gamma'<\gamma'
\end{align*}
where the inequality is due to the assumption
(\ref{eq:fixed-cost}). For this $\gamma^{'}$, we have
$V^{\gamma^{'}}_a(a)<\gamma^{'}$.

Since $\gamma \rightarrow V^{\gamma}_a$ is continuous, which
follows from the fact that this function is convex, and
increasing, $V^0_a>0$ and $V^{\gamma^{'}}_a(a)<\gamma^{'}$ implies
that $\gamma \rightarrow V^{\gamma}_a$ crosses the line $\gamma
\rightarrow \gamma$.

\end{proof}

\begin{lemma}\label{lem:majorant-forex}
Assume that
\begin{equation}\label{eq:defn-mxa}
r(x,a) \quad \text{is lower semi-continuous.}
\end{equation}
 Let us define $R^{\gamma}(\cdot;
a)\triangleq\frac{r^{\gamma}(F^{-1}(\cdot),a)}{\varphi(F^{-1}(\cdot))}$
where
\begin{equation}\label{eq:appendix-tilde-r}
r^{\gamma}(x, a)\triangleq \ME[e^{-\alpha\Delta}(\bar{K}(X_\Delta,
a)+\gamma)].
\end{equation}
Then (\ref{eq:Wa}) is the smallest non-negative concave majorant
of $R^{\gamma}$ that passes through  $(F(c+),l_c)$.
\end{lemma}

\begin{proof}
See for e.g. Dynkin \cite{dynkin} and Dayanik and Karatzas
\cite{DK2003}.
\end{proof}
\begin{lemma}\label{lemma-excessive}
Assume that (\ref{eq:positivity}) and (\ref{eq:defn-mxa}) hold.
Then $u_a/\varphi$ is $F-$concave, i.e.,
$\alpha-$excessive.\footnote{A function $f$ is called
$\alpha$-excessive function of $X_0$ if for any stopping time
$\tau$ of the natural filtration of $X^0$ and $x \in (c,d)$, $f(x)
\geq  \E^{x}\left[ e^{-\alpha \tau} f(X^0_{\tau})\right]$, see for
e.g. \cite{salminen} and \cite{dynkin} for more details.}
\end{lemma}
\begin{proof}
This follows from Lemmas~\ref{lem:existence-lemma} and
\ref{lem:majorant-forex}. For the equivalence of
$\alpha$-excessivity and $F-$concavity see e.g. Theorem 12.4 in
\cite{dynkin} and also \cite{DK2003}. This fact can be observed
from (\ref{eq:devH}).
\end{proof}
\begin{lemma}\label{eq:u-sub-u-sup-ineq}
Assume that (\ref{eq:positivity}) and (\ref{eq:defn-mxa}) hold.
Then
\begin{equation}\label{eq:u-sup-sub}
u^a(x) \leq u_a(x), \quad x \in (c,d).
\end{equation}
\end{lemma}
\begin{proof}
It follows from Lemma~\ref{lemma-excessive} that $u_a$ is
$\alpha$-excessive. Also, observe from
(\ref{eq:inner-optimization}) that
\begin{equation}\label{eq:ua-at0}
u_a(x) \geq r(x;a)+e^{-\alpha \Delta} u_a(a),
\end{equation}
where $r(x,a)$ is as in (\ref{eq:defn-mxa}). Let $\nu=\{T_1, T_2,
...,T_i,...;\xi_1, \xi_2,...,\xi_i,...\}$ be an admissible control
and let $T_0=0$. Without loss of generality we will assume that
$r(b;a)>0$, because otherwise the corresponding strategy will have
a lower value function $J^{\nu}(x)$ associated to it. Since $u_a$
is $\alpha-$ excessive,
\begin{equation}
\begin{split}
u_a(x) \geq \E^x\left[e^{-\alpha T_1} u_a(X_{T_1})\right], \quad
\text{and} \quad \E^x\left[e^{-\alpha
(T_i+\Delta)}u_{a}(X_{(T_i+\Delta)})\right]-\E^x \left[e^{-\alpha
T_{i+1}}u_a(X_{T_{i+1}})\right] \geq 0,
\end{split}
\end{equation}
for all $i=1,...,N-1$. Then
\begin{equation}
\begin{split}
u_a(x) &\geq \E^x\left[e^{-\alpha T_1}
u_a(X_{T_1})\right]+\sum_{i=1}^{N-1}\E^x \left[e^{-\alpha
T_{i+1}}u_a(X_{T_{i+1}})\right]-\E^x\left[e^{-\alpha
(T_i+\Delta)}u_{a}(X_{(T_i+\Delta)})\right]
\\& =\E^x \left[e^{-\alpha
T_{N}}u_a(X_{T_{N}})\right]+\sum_{i=1}^{N-1}\E^x \left[e^{-\alpha
T_{i}}u_a(X_{T_{i}})\right]-\E^x\left[e^{-\alpha
(T_i+\Delta)}u_{a}(X_{(T_i+\Delta)})\right]
\\&\geq \sum_{i=1}^{N-1}\E^{x}\left[e^{-\alpha T_i}r(X_{T_i},a)\right],
\end{split}
\end{equation}
in which the inequality follows from (\ref{eq:ua-at0}) and the
fact that $u_a$ is non-negative. Now, using the monotone
convergence theorem
\begin{equation}
\begin{split}
u_a(x) &\geq \E^{x}\left[\sum_{i=1}^{\infty}e^{-\alpha T_i}
r(X_{T_i},a)\right] =\E^x\left[\sum_{i=1}^{\infty}
e^{-\alpha(T_i+\Delta)}\E^{X_{T_i}}\left[\bar{K}(X_{\Delta},a)\right]\right]
\\&=\E^{x}\left[\sum_{i=1}^{\infty}e^{-\alpha(T_i+\Delta)}\E^{X_{T_i}}\left[K(X_{\Delta},a)-g(X_{\Delta})+g(a)\right]
\right]
\\&=\E^{x}\left[\sum_{i=1}^{\infty}e^{-\alpha(T_i+\Delta)}K(X_{(T_i+\Delta)-},X_{T_i+\Delta})\right]
+\E^{x}\left[\sum_{i=1}^{\infty}e^{-\alpha(T_i+\Delta)}(-g(X_{(T_i+\Delta)-})+g(X_{(T_i+\Delta)}))\right]
\\&=\E^{x}\left[\sum_{i=1}^{\infty}e^{-\alpha(T_i+\Delta)}K(X_{(T_i+\Delta)-},X_{T_i+\Delta})\right]
+\E^{x}\left[\int_0^{\infty}e^{-\alpha s}f(X_s)ds\right]-g(x)
=u^{a,b}(x).
\end{split}
\end{equation}
The third inequality follows from Remark~\ref{rem:not-Markov}).
The fourth inequality can be derived from
(\ref{eq:integral-of-f}). The last equality follows from
(\ref{eq:defn-u-J-g}). Now taking to supremum over $b$, we obtain
(\ref{eq:u-sup-sub}).
\end{proof}

\begin{lemma}\label{lem:u-sub-a-eq-u-sup-a}
Assume that (\ref{eq:positivity}) and (\ref{eq:defn-mxa}) hold and
that the function $x\rightarrow R(x;a)$ defined in
(\ref{eq:defn-r-forex}) is concave and increasing on $(a',d)$ for
some $a' \in (a,d)$ and that
\begin{equation}\label{eq:infty}
\lim_{x \rightarrow F(d)}R(x;a)=\infty.
\end{equation}
Then $u_a(x)=u^{a, b^*}(x)$ for a unique $b^*\in (c,d)$. Hence
from Lemma~\ref{eq:u-sub-u-sup-ineq} it follows that
$u_a(x)=u^{a}(x)=u^{a ,b^*}(x)$, $x \in (c,d)$.
\end{lemma}
\begin{proof}
Since $R$ is concave, $R^{\gamma}$ in (\ref{eq:appendix-tilde-r})
is also concave on $(a', d)$. The assumption in (\ref{eq:infty})
implies that the smallest concave majorant $W^{\gamma}_a$ in
(\ref{eq:Wa}) is linear on $(F(c),F(b^\gamma))$ for a
\emph{unique} $b^{\gamma} \in (c,d)$ and is tangential to
$R^{\gamma}(\cdot,a)$ at  $F(b^{\gamma})$ and coincides with
$R^{\gamma}(\cdot,a)$ on $[F(b^{\gamma}),F(d))$. Together with
Lemma~\ref{lem:existence-lemma} this implies that there exists a
unique $\gamma^*$ such that
 equations (\ref{eq:rem1}) and
(\ref{eq:rem2}) are satisfied when $W$ is replaced by
$W_a^{\gamma^*}$ and $b$ is replaced by $b^{\gamma^*}$. Note that
$W_a^{\gamma^*}$ corresponds to a strategy $(a,b^{\gamma^*})$.
That is, if we start with $u^{a,b^{\gamma^*}}$ and transform it
via (\ref{eq:W-u-trnasfrom}) we get $W_a^{\gamma^*}$. On the other
hand, using (\ref{eq:Wa}) with by substituting $\gamma=\gamma^*$
we have that $u_a(x)=\varphi(x)W_a^{\gamma^*}(F(x))$, $x \in
(c,d)$. This let's us conclude that $u^{a,b^{\gamma^*}}=u_a(x)$,
$x \in (c,d)$. We see that the unique $b^*$ in the claim of the
proposition is $b^{\gamma^*}$.
\end{proof}

\begin{proposition}\label{prop:existence-uniquenes-forex}
Assume that the hypotheses of Lemma~\ref{lem:u-sub-a-eq-u-sup-a}
are satisfied. Then there exists a unique solution to
(\ref{eq:non-linear-ipmlicit-FOREX}). If $b^{*}$ is the unique
solution of (\ref{eq:non-linear-ipmlicit-FOREX}), then
$u^{a}(x)=u^{a,b^*}(x)$.
\end{proposition}
\begin{proof}
In the proof of Lemma~\ref{lem:u-sub-a-eq-u-sup-a}, we have seen
that there exists a unique $b^*$ such that (\ref{eq:rem1}) and
(\ref{eq:rem2}) are satisfied. Using
Remark~\ref{rem:smooth-fit-forex}, we conclude that $b^*$ is the
unique solution of (\ref{eq:non-linear-ipmlicit-FOREX}).
\end{proof}

Note that when the assumptions of Proposition
\ref{prop:existence-uniquenes-forex} hold, the optimal threshold
strategy is described by a single open interval in the state space
of the controlled process. The conditions for the existence and
uniqueness of the optimal interval are specified, essentially by
the conditions on total reward function $\bar{K}(x, y)$ associated
with one intervention from $x$ to $y$ (see (\ref{eq:u1-under}),
(\ref{eq:R-forex}) ) and drift and volatility of the underlying
diffusion as the function $F$ depends on them that appears in
(\ref{eq:R-forex}) depends on them.

\underline{Third stage}: Now, we let $a\in (c,d)$ vary and choose
$a^*$ that maximizes $\rho(a)$ and also find $b^*=b(a^*)$.
\, Finally, we obtain the value function given in (\ref{problem1})
by $v(x)=u(x)+g(x)$.

\subsection{Example: Optimal Exchange Rate Intervention When There is
Delay}\label{sec:example-one-sided} To illustrate the procedure of
solving impulse control problems with delay, we take an example
from Mundaca and {\O}ksendal \cite{mundaca-oksendal} (also see
{\O}ksendal \cite{O1999}) that considers the following foreign
exchange rate intervention problem:
\begin{equation}\label{oksprob}
   J^\nu_D(x) \triangleq \ME\left[\int_0^\infty e^{-\alpha s}X_s^2 ds +
   \sum_{i}^{\infty}e^{-\alpha (T_i+\Delta)}(c+\lambda|\xi_i|)\right]
\end{equation}
where $X^0_t=x+B_t$, in which $B$ is a standard Brownian motion.
Here, the superscript 0 is to indicate that the dynamics in
consideration are of the uncontrolled state variable. In
(\ref{oksprob}), $c>0$ and $\lambda\geq 0$ are constants
representing the cost of making an intervention. The problem
without delays are solved by \cite{O1999} through
quasi-variational inequalities and by \cite{DE2005} using a direct
characterization of the value function. In this problem, the
Brownian motion represents the exchange rate of currency and the
impulse control represents the interventions the central bank
makes in order to keep the exchange rate in a given target window.
At time $T_i$, such that $X_{T_i-}=b$, the central bank makes a
commitment to reduce the exchange rate from $b$ to $a<b$, which is
implemented $\Delta$ units of time later. During the time interval
$(T_i, T_i+\Delta]$ the central bank does not make any other
interventions. $\Delta$ units later if the exchange rate is still
greater than $a$, then the central bank reduces the exchange rate
from $X_{(T_i+\Delta)-}$ to $a$ and pays a cost of $c+\lambda
(X_{(T_i+\Delta)-}-a)$. On the other hand, if $\Delta$ units of
time later if the exchange rate is less than $a$, the central bank
chooses increases the exchange rate to $a$ at a cost of $c+\lambda
(a-X_{(T_i+\Delta)-})$. This is a one-sided impulse control
problem, in the sense that a control is triggered only if $X_t>b$
and there has not been any previous action in the interval
$(t-\Delta,t)$.

 The problem is to minimize the expected total
discounted cost over all threshold strategies.
\begin{equation}\label{eq:minimization}
v_D(x) \triangleq \inf_{\nu} J_D^\nu(x).
\end{equation}
A similar version of this problem is analyzed by
{\O}ksendal and Sulem \cite{oksendal-delay-impulse}, in which they
take the controls $\xi_i \in \mathcal{F}_{T_i}$ for all $i$. (This
introduces path dependence since the value of $X_{T_i+\Delta}$ is
partially determined by $\mathcal{F}_{T_i}$.)

Instead of
solving a minimization problem of (\ref{eq:minimization}), we will
solve
\begin{equation*}
v(x)=\sup_{\nu} \ME\left[\int_0^\infty e^{-\alpha s}(-X_s^2) ds -
   \sum_{i}^{\infty}e^{-\alpha (T_i+\Delta)}(c+\lambda|\xi_i|)\right].
\end{equation*}
and recover the value function by $v_D(x)=-v(x)$. (Here, the supremum is taken over all the threshold strategies.) The continuous
cost rate is $f(x)=-x^2$ and the intervention cost is $K(x,
y)=-c-\lambda|x-y|$ in our terminology. By solving the equation
$(\mathcal{A}-\alpha)v(x)=\frac{1}{2}v^{''}(x)-\alpha v(x)=0$, we
find that $\psi(x)=e^{x\sqrt{2\alpha}}$ and
$\varphi(x)=e^{-x\sqrt{2\alpha}}$.  Hence
$F(x)=e^{2x\sqrt{2\alpha}}$ and $F^{-1}(x)=\frac{\log
x}{2\sqrt{2\alpha}}$.   
Using Fubini's theorem we can calculate $g(x)$ explicitly as:
\begin{align*}
g(x)=-\ME\int_0^\infty e^{-\alpha
s}(x+B_s)^2ds=-\left(\frac{x^2}{\alpha}+\frac{1}{\alpha^2}\right).
\end{align*}

\noindent We shall follow the procedure described in the last
section: Let us fix $a>0$ and consider
\begin{align}\label{eq:r-forex-example}
&r(x, a)=\ME[e^{-\alpha\Delta}\bar{K}(X_\Delta, a)]
=\ME\left[e^{-\alpha\Delta}\Big(-c-\lambda|X_\Delta-a|+g(a)-g(X_\Delta)\Big)\right]\\
&=\ME\left[e^{-\alpha\Delta}\left(-c-\lambda|x+B_\Delta-a|
-\left(\frac{a^2}{\alpha}+\frac{1}{\alpha^2}\right)+\left(\frac{(x+B_\Delta)^2}{\alpha}
+\frac{1}{\alpha^2}\right)\right)\right] \nonumber \\
&=e^{-\alpha\Delta}\left(-c-\lambda \left(2 \Delta
\exp\left(-\frac{(a-x)^2}{4 \Delta^2}\right)+ (a-x)\left(-1+2
N\left(\frac{a-x}{\Delta}\right) \right)\right)
+\frac{x^2-a^2+\Delta}{\alpha}\right) \nonumber .
\end{align}
The left boundary $-\infty$ is natural for a Brownian motion and,
for any $a>0$,
\begin{equation*}\label{}
    l_{-\infty}=\limsup_{x\downarrow-\infty}\frac{r(x, a)^+}{\varphi(x)}=0.
\end{equation*}
It follows that $R(y)$ passes through $(F(-\infty),
l_{-\infty})=(0, 0)$. (See Dayanik and Karatzas\cite{DK2003}
Proposition 5.12.)

\begin{proposition}\label{prop:forex-example}
For the function $r$ in (\ref{eq:r-forex-example}), there exists a
unique
solution to (\ref{eq:non-linear-ipmlicit-FOREX}) for a fixed $a$.
\end{proposition}
\begin{proof}
See Appendix.
\end{proof}

Using the algorithm we described in Section~\ref{sec:algorithm-threshold}
we find the optimal $(a^*, b^*, \rho^*)$. Going
back to the original space we get
\begin{equation*}
V(x)=\sup_{a, b\in \R}
u(x)=\varphi(x)W^*(F(x))=\varphi(x)(\beta^*)F(x)=\rho^*
e^{x\sqrt{2\alpha}}.
\end{equation*}
on $x\in(-\infty, b^*]$. To get $v(x)=\sup_{\nu}J^\nu(x)$, we add
back $g(x)$,
\begin{equation*}
v(x)=V(x)+g(x)=\rho^*e^{x\sqrt{2\alpha}}-\left(\frac{x^2}{\alpha}+\frac{1}{\alpha^2}\right).
\end{equation*}
Finally, flipping the sign we obtain the optimal cost function as
\begin{eqnarray}\label{eq:for-ex-value-func}
v_D(x)&=& \begin{cases}
                    \hat{v}_o(x)\triangleq\left(\frac{x^2}{\alpha}+\frac{1}{\alpha^2}\right)-\rho^*e^{x\sqrt{2\alpha}}, &0\leq x \leq b^*, \\
                   - e^{-\alpha\Delta}\rho^*e^{a^*\sqrt{2\alpha}}-r(x; a^*)+
                    \frac{x^2}{\alpha}+\frac{1}{\alpha^2}, &b^*\leq x.
        \end{cases}
\end{eqnarray}
Figure 1 is obtained when the parameters are chosen to be $(c,
\lambda, \alpha, \Delta)=(150, 50, 0.2, 1.0)$. We found the
solution triplet to be $(a^*, b^*, \rho^*)=(5.066, 12.1756,
0.042423)$.
 The optimal cost function without delay, for the
same parameters, has the solution triplet $(a_0, b_0,
\rho_0)=(5.07723, 12.2611, 0.0492262)$. The continuation region
shifts to the left with delay (it shrinks from $(-\infty,
12.2611)$ to $(-\infty, 12.1756)$), and the central bank acts more
aggressively when it encounters delays (see Figure 1 - (c)).
\begin{figure}[h]\label{figure:1}
\begin{center}
\begin{minipage}{0.45\textwidth}
\centering \includegraphics[scale=0.75]{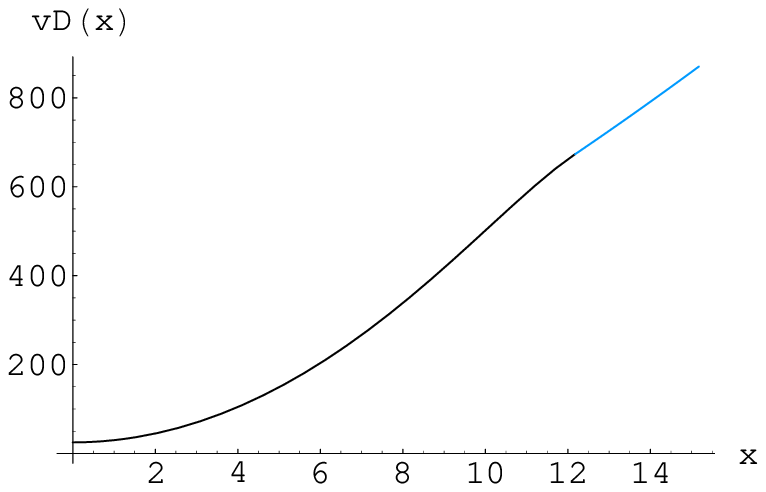} \\
(a)
\end{minipage}
\begin{minipage}{0.45\textwidth}
\centering \includegraphics[scale=0.75]{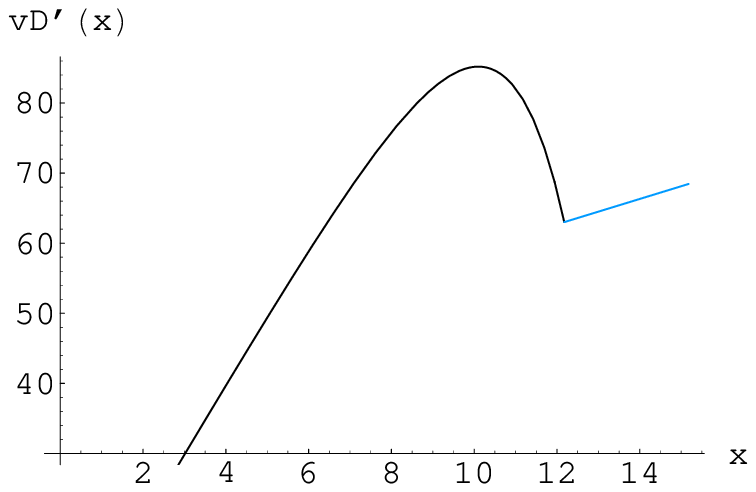} \\
(b)
\end{minipage}
\begin{minipage}{0.45\textwidth}
\centering \includegraphics[scale=0.75]{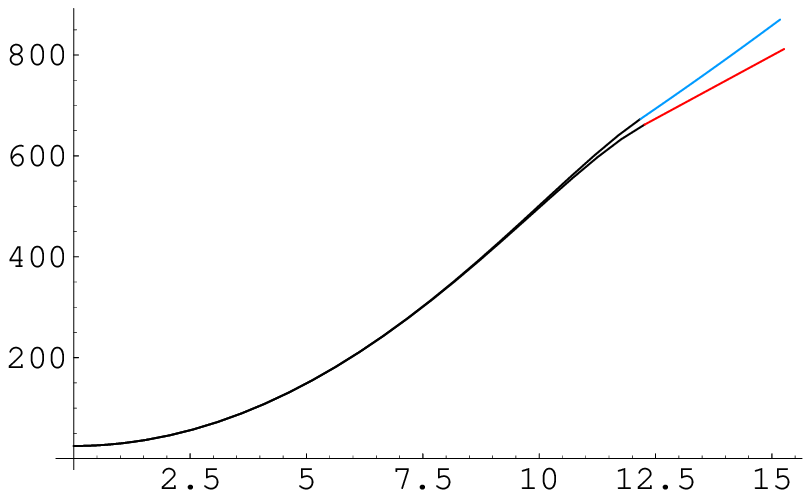} \\
(c)
\end{minipage}
\begin{minipage}{0.45\textwidth}
\centering \includegraphics[scale=0.75]{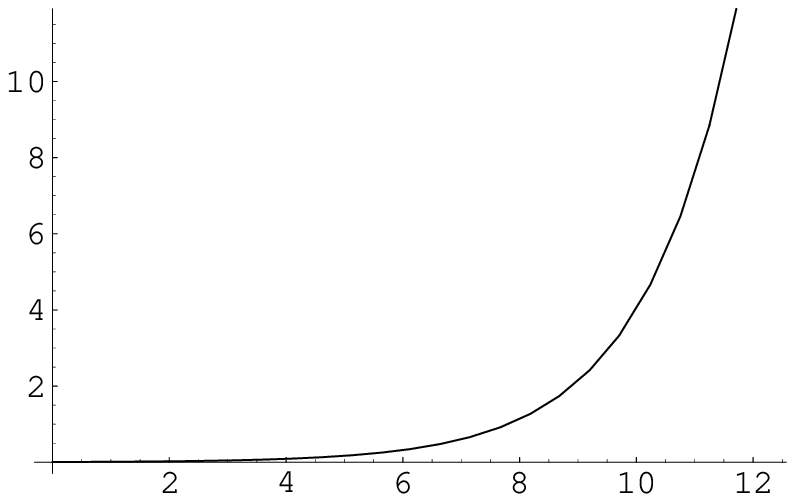} \\
(d)
\end{minipage}
\caption{\small (a) The optimal cost function $v_D(x)$. The dotted
line and the solid line fit each other continuously at
$b^*=12.1756$. (b) The derivative of $v_D(x)$, showing that the
smooth-fit principle holds at $b^*$. (c) Comparison of $v_D(x)$
with the cost function without delay $v_0(x)$. Note that $v_D$
majorizes $v_0$. (d) Plot of the difference of $v_D(x)-v_0(x)$.}
\end{center}
\end{figure}

\section{Firing Costs and Labor Demand: Optimal Band Strategies}\label{sec:firing} In this section, we will improve on the
techniques of the previous section in order to study an impulse
control corresponding to band policies when there are
implementation delays. In particular, we will concentrate our
attention on a specific example, which is of practical interest.
We will find optimal hiring and firing decisions of a firm that
faces stochastic demand and has to conform to regulatory delays
when it is firing employees.

Recently, General Motors Corporation (GM) has decided to lay off
25,000 of its work force to cut back on its production and administrative costs. However ``GM's
UAW (United Auto Workers) contract essentially forces it to pay
union employees during the life of the contract even if hourly
workers are laid off and their plants are closed. But those
protections only run
 through September 2007, when the current four-year pact with
the union ends. GM spokesman Ed Snyder said the automaker has yet
to reach any agreement with the UAW yet on the nature or the
manner of the work force reduction." \footnote{Source:  June 7,
2005 CNN Money, ``GM to cut 25,000 jobs" by Chris Isidore,
http://money.cnn.com/2005/06/07/ news/fortune500/gm\_closings/}
This is a typical example of a firing cost and implementation
delay a corporation faces when the workers are unionized. Another
example of firing delay caused by government regulations in Europe
(see e.g. \cite{bb}).

Bentolila and Bertola \cite{bb} address the issue of costly hiring
and firing and its effects on unemployment rate in Europe using
singular stochastic control. Here, we are solving an impulse
control problem since we are also taking fixed cost of labor
adjustments into account. But our main purpose is to measure the
effects on firing delay in decisions of firms. As we shall see, it
turns out that the controlled state variable is not Markov,
therefore we will focus our attention completely on the \emph{band
policies} (which we will define shortly) rather than trying to
find the best impulse control policy. Our method of solving
impulse control problem differs from its counterparts that use
quasi-variational inequalities since we give a direct
characterization of the value function as a linear function in the
\emph{continuation} region without having to guess the form of the
solution and without having to prove that the conjectured solution
satisfies conditions of a verification lemma.

\subsection{Problem setup}

As in \cite{bb}\footnote{The set up of Bentolila and Bertola
\cite{bb} was brought to our attention by Keppo and Maull. In the
INFORMS Annual Meeting in 2004, Keppo and Maull presented their
partial results on the hiring and firing decisions of firms which
they obtained by solving quasi-variational inequalities.}, we will
consider a firm with a linear production technology. In particular
the quantity sold is $Q_t=A L_t$, $A \in \mathbb{R}_+$, in which
$L_t$ is the labor at time $t$. The selling price at time $t$,
$P_t$, of the product  is determined from
\begin{equation}\label{eq:quantity}
Q_t= Z_t P_{t}^{\frac{1}{\mu-1}}, \quad \mu \in (0,1)
\end{equation}
in which $Z_t$ indexes the position of the direct demand curve
whose dynamics follow
\begin{equation}\label{eq:demand}
dZ_t= Z_t b dt + Z_t \sigma_t dW_t
\end{equation}
with a constant $b\in \R_+$.  In equation (\ref{eq:quantity}) the
quantity $1-\mu$ is the firm's monopoly power. Let us denote the
filtration generated by
 the demand process $Z$ by $\mathcal{F} \triangleq
(\mathcal{F}_{t})_{t \geq 0}$. We will make the following
assumption to guarantee that (\ref{eq:demand}) has a unique strong
solution. We assume that $\sigma$ is bounded and adapted to the
filtration of the Brownian motion $W$.

 In our framework, if the firm
produces excess products because of the excess labor, the products
produced are still all sold but at a cheaper price. The firm pays
a wage, $w$, to its workers, therefore the net rate of profit that
the firm makes at time $t$ is given by
\[
Q_t P_t-wL_t = Z_t^{1-\mu} (A L_t)^{\mu} - w L_t.
\]
When the workers quit voluntarily, the firm bears no firing costs
and we assume that the workers quit at rate $\delta$, that is,
without any intervention from the management the labor force
follows the dynamics
\begin{equation}\label{eq:L0}
dL^0_t=-\delta L^0_t dt.
\end{equation}
Here, as in the previous section, the superscript 0 indicates that
there are no controls applied. The firm makes commitments to
change its labor force at times $\{S_i\}_{i \in \mathbb{N}}$ and
$\{T_i\}_{i \in \mathbb{N}}$. At time $S_i$ the firm makes a
commitment to increase its labor force (which is immediately
implemented), and at time $T_i$ it makes a commitment to decrease
its labor force, which is implemented $\Delta$ units of time
later. During the time interval $(T_i,T_i+\Delta]$ the firm makes no
commitments to change its labor force. Note that although at time
$T_i$ the firm decided to decrease its labor force, the labor
force itself might move to very low levels following the dynamics
(\ref{eq:L0}), therefore at time $T_i+\Delta$ the firm may end up
hiring to move keep the production level up. However, if the labor force level
is still very high at time $(T_i+\Delta)-$, then the firm ends up
firing. Here, $\Delta$ represents the regulatory delays a firm
faces when it is cutting off its work force.

The labor adjustments come at a cost: At time $S_i$ the firm
increases the labor by $\zeta_i (\geq 0) \in \mathcal{F}_{S_i}$
(Here, for the sake of brevity we are taking the $\sigma$-algebras
as a collection of mappings.) to $L_{S_i-}+\zeta_i$, then the
associated cost is
\[
c_1 \zeta_i + c_2 L_{S_i-}.
\]
At time $T_i$, the  firm makes a commitment to decrease the labor
at time $T_{i}+\Delta$. If it ends up decreasing the labor force
by $\eta_i (\geq 0) \in \mathcal{F}_{T_i+\Delta}$ to
$L_{T_{i}+\Delta}= L_{(T_i+\Delta)-}-\eta_i$, then the associated
cost is quantified as
\[
c_3 \eta_i + c_4 L_{(T_i+ \Delta) -},
\]
which depends on the amount of labor force to be fired and the
level of the total labor force as well.  The latter component of
costs is based on the following observations:  When a corporation
decides who to be fired or which division to be restructured,
administrative costs will become larger in proportion to the size
of the total labor force since the firm's operations are closely
knitted among various divisions.

On the other hand as we discussed above if the labor force itself
moves to very low levels itself during the $\Delta$ units of time,
 at time $T_i+\Delta$ the firm may end up hiring (in this
case $\eta_i \leq 0$) to keep the production up at the cost of
\[
c_1 |\eta_i|+c_2 L_{(T_i+ \Delta) -}
\]
for some positive constants $c_1, c_2, c_3, c_4$ and $\Delta \geq
0$. This cost becomes negligible as $\Delta$ becomes small because
in that case the work force does not change much by itself. So the
controls of the firm are of the form
\[
\nu=( S_1, S_2, \cdots; \zeta_1, \zeta_2, \cdots; T_1, T_2,
\cdots; \eta_1,\eta_2, \cdots),
\]
where $0 \leq S_1 < S_2< \cdots$ and $0 \leq T_1 < T_2 < \cdots$
are two increasing sequences of stopping times of the filtration
$\mathcal{F}$. $T_{i+1}-T_i \geq \Delta$ and for any $i$ there
exists no $j$ such that $T_i \leq S_j \leq T_{i+\Delta}$. The
magnitudes of the impulses satisfy $\zeta_i (\geq 0) \in
\mathcal{F}_{S_i}$ and $\eta_i (\in \mathbb{R}) \in
\mathcal{F}_{T_i+\Delta}$ for all $i$. We call these type of
controls \emph{admissible} and we will denote the set of all
admissible controls by $\mathcal{V}$. To each control $\nu \in
\mathcal{A}$ we associate a profit function of the form
\begin{equation}
\begin{split}
 J^{\nu}(z,l)& \triangleq \mathbb{E}\bigg[\int_{0}^{\infty} e^{-rt}
\left(Z_{t}^{1-\mu} (A L_t)^{\mu} -w L_t \right)dt  - \sum_{i}
e^{- r S_i} \left(c_1 \zeta_i + c_2 L_{S_i-}\right)
\\ & - \sum_{j} e^{-r(T_j + \Delta)} \left(\left(c_3 \eta_j + c_4
L_{(T_j+\Delta)-}\right) 1_{\{\eta_j>0\}}+ \left(c_1 \eta_j + c_2
L_{(T_j+\Delta)-}\right) 1_{\{\eta_j<0\}}\right) \bigg],
\end{split}
\end{equation}
which incorporates the profit and cost structure we described so
far. Here $r>b$ is a subjective rate of return that the firm uses
to discount its future profits. In fact if $r<b$, then taking no
action is optimal as we will point out below. Under the measure
$\mathbb{P}$,  we have that $L_0=l$ and $Z_0=z$ almost surely.

The objective of the company is then to maximize its profits by
choosing the best possible strategy $\nu^*$ such that
\begin{equation}\label{eq:value}
v(z,l) \triangleq \sup_{v \in \mathcal{V}} J^{\nu} (z,l)
=J^{\nu^*}(z,l),
\end{equation}
if the optimal strategy $\nu^*$ exists. Hereafter, we will refer
to $v$ as the value function.

It looks as if the control problem defined in (\ref{eq:value})
involves two state variables, namely the demand $Z$ and the labor
force $L$. Recall that we have no control over the demand $Z$ but
we can control the labor force $L$ by making hires and fires. But
the only source of randomness is the demand process. In the sequel
we will show that the optimal control problem (\ref{eq:value})
involves only one state variable. On denoting $\xi_t \triangleq
L_t/Z_t$, $t \geq0$ and the absolute changes in labor per unit of
demand by $\beta_i \triangleq \zeta_i/Z_{S_i}$ and $\alpha_i
\triangleq \eta_i/ Z_{T_i+\Delta} \in \mathcal{F}_{T_i+\Delta}$,
we can write the the profit function $J^{\nu}$ as
\begin{equation}
\begin{split}
&J^{\nu}(z,l)=\mathbb{E}\bigg[\int_0^{\infty}e^{-rt}Z_t\left((A
\xi_t)^{\mu}- w \xi_t \right)dt- \sum_{i} e^{-r S_i} Z_{S_i -}(c_1
\beta_i+ c_2 \xi_{S_i -}) \\ & \sum_{j} e^{-r(T_j + \Delta)}\left(
(c_{3} Z_{T_j}\, \alpha_j + c_4 Z_{T_j +\Delta}\,\xi_{(T_j +
\Delta)-})1_{\{\alpha_j>0\}}+(c_{1} Z_{T_j}\, \alpha_j + c_2
Z_{T_j +\Delta}\,\xi_{(T_j + \Delta)-})1_{\{\alpha_j<0\}}\right)
\bigg].
\end{split}
\end{equation}
Let us introduce a new probability measure $\mathbb{P}_0$ by
\begin{equation}
\frac{d\mathbb{P}_0}{d\mathbb{P}} \bigg|_{\mathcal{F}_t}=
\tilde{Z}_t, \quad \text{where}\quad  \tilde{Z}_t=
\exp\left(\int_0^{t} \sigma_s dW_s - \frac{1}{2} \int_0^t
\sigma_s^2ds \right)
\end{equation}
for every $0 \leq t< \infty$. Using the representation of the
profit function $J^{\nu}$, we can write it as
\begin{equation}\label{eq:Jnu-eq-I}
J^{\nu}(z,l)=z I^{\nu}\left(\frac{z}{l}\right),
\end{equation}
in which
\begin{equation}
\begin{split}
I^{\nu}(\xi) & \triangleq
\mathbb{E}^{\xi}_0\bigg[\int_{0}^{\infty}e^{(b-r)t} z \left((A
\xi_t)^{\mu} - w \xi_t \right)dt - \sum_i e^{(b-r)S_i} (c_1
\beta_i+ c_2 \xi_{S_i -}) \\&- \sum_j e^{(b-r)(T_j + \Delta)}
 \left((c_{3} \alpha_j+ c_4 \xi_{(T_j + \Delta)-})1_{\{\alpha_j>0\}}
 +(c_{1} \alpha_j+ c_2 \xi_{(T_j + \Delta)-})1_{\{\alpha_j<0\}} \right)\bigg],
\end{split}
\end{equation}
where $\mathbb{E}^{\xi}$ is the expectation under $\mathbb{P}_0$
given that $\xi_0=\xi$. Here, with slight abuse of notation, on
the right-hand-side of (\ref{eq:Jnu-eq-I}), we denoted
\[
\nu=( S_1, S_2, \cdots; \beta_1, \beta_2, \cdots; T_1, T_2,
\cdots; \alpha_1,\alpha_2, \cdots),
\]
is a control that is applied to the process $\xi$. The controls
here are such that $\beta_i (\geq 0) \in \mathcal{F}_{S_i}$ and
$\alpha_i (\in \mathbb{R}) \in \mathcal{F}_{T_i+\Delta}$.  Again
as before $\{S_n\}_{n \in \mathbb{N}}$ and $\{T_n\}_{n \in
\mathbb{N}}$ are two increasing sequence of stopping times. We
also assume that $T_{i+1}-T_i \geq \Delta \geq 0$ and that for any
$i$ there exists no $j$ such that $T_i \leq S_j \leq
T_{i+\Delta}$. With another slight abuse of notation we will
denote the admissible set of controls we described here also by
$\mathcal{V}$. As a result of the developments in the last part of
this section we see that the process $L_t/Z_t$ is the sufficient
statistic of the problem in (\ref{eq:value}). In fact we can write
the value function as
\begin{equation}\label{eq:Y}
v(z,l)=z Y\left(\frac{z}{l}\right),\,\,\,  \text{where} \quad
Y(\xi) \triangleq \sup_{\nu \in \mathcal{V}}I^{\nu}(\xi).
\end{equation}
Under the measure $\mathbb{P}^{\xi}_0$ the dynamics of the
process, $\xi_t$ when there are no impulses applied follows
\begin{equation}\label{eq:xi0t}
\xi^0_t=\xi \exp\left(-(b+\delta)t -\int_{0}^{t} \sigma_s dB_s-
\frac{1}{2} \int_0^{t} \sigma^2_s ds \right),
\end{equation}
where $B$ is a Wiener process under measure $\mathbb{P}_0$. Here,
as before, the superscript 0 indicates that there are no
controls/impulses applied.

\subsection{Solution}

Although the controlled process $\xi$ is not a Markov process,
because depending on whether the process reaches a point during
the interval $(T_i, T_i+\Delta)$ or not, that point has different
roles. That is, how the process reaches to a
particular point (path information) affects how the process will
continue from this point. However, the process regenerates at
times $\{T_i+\Delta\}_{i \in \mathbb{N}}$ and the value of the
process at time $T \in (T_i,T_{i+\Delta})$, $X_T$, depends on the
information up to $T_i$, $\mathcal{F}_{T_i}$, only through the
value of the process at time $T_i$, $X_{T_i}$. Therefore, as we did in Section
\ref{sec:characterization}, assuming there is no
history prior to time 0, i.e. $\mathcal{F}_{0}$ is a trivial
sigma-algebra, we can develop
\begin{equation}
\begin{split}
I^{\nu}(\xi)&=\mathbb{E}^{\xi}_0\bigg[1_{\{T_1<S_1\}}e^{(b-r)(T_1+\Delta)}\left(C_1(\xi_{(T_1+\Delta)-
},\xi_{T_1+\Delta})-g(\xi_{(T_1+\Delta) -})+ I^{\nu}
(\xi_{T_1+\Delta})\right)
\\ & + 1_{\{T_1>S_1\}} e^{(b-r)S_1}\left(C_2(\xi_{S_1
- },\xi_{S_1})-g(\xi_{S_1-})+ I^{\nu} (\xi_{S_1})\right)\bigg],
\end{split}
\end{equation}
\begin{equation}
\begin{split}
\text{where}\quad & C_2(x,y) \triangleq -c_1(y-x)1_{\{y>x\}}- c_2
x, \quad \text{and},
\\ & C_1(x,y)\triangleq - (c_3 (x-y)+c_4 x) 1_{\{x>y\}}+ C_2(x,y)
1_{\{y>x\}}
\end{split}
\end{equation}
\begin{equation}\label{eq:g-labor}
g(\xi)\triangleq \mathbb{E}_0\left[\int_{0}^{\infty}e^{(b-r)t}
\left((A \xi^0_t)^{\mu} - w \xi^0_t \right)dt\right].
\end{equation}
On denoting $u(\xi) \triangleq I^{\tilde{\nu}}(\xi)-g(\xi)$, we
can write
\begin{equation}\label{eq:u}
\begin{split}
u(\xi)&=\mathbb{E}^{\xi}_0\left[1_{\{T_1<S_1\}}e^{(b-r)(T_1+\Delta)}
\left(\bar{C}_1(\xi_{(T_1+\Delta) - },\xi_{T_1+\Delta})+
u(\xi_{T_1+\Delta}) \right) \right]
\\ &+ \mathbb{E}_0^{\xi}\left[1_{\{T_1>S_1\}}e^{(b-r)S_1}
\left(\bar{C}_2(\xi_{S_1- },\xi_{S_1})+ u(\xi_{S_1}) \right)
\right],
\end{split}
\end{equation}
in which
\begin{equation}
\bar{C}_1(x,y) \triangleq C_1(x,y)-g(x)+g(y) \quad \text{and}\quad
\bar{C}_2(x,y) \triangleq C_2(x,y)-g(x)+g(y).
\end{equation}

In the rest of this section, we will analyze the following double
sided threshold strategy (\emph{band policy}) of the following
form: 1) Whenever the marginal revenue product of labor hits level
$d$, the firm makes a commitment to bring the
marginal revenue product of labor to $c<d$. This may be achieved
by firing employees if marginal revenue product of labor is still
greater than $c$ after the delay. However, it is possible that
after the delay the marginal revenue product of labor will be less than
$c$. In this case, the firm makes hires.  2) Whenever the marginal revenue
product of labor hits level $p$ the firm increases it to $q>p$ (by
hiring new employees). We will characterize the value
function corresponding to an arbitrary band policy.

For a band policy we described above $S_{1}=\tau_p$ and
$T_{1}=\tau_d$, and
\[
\xi_{T_1+\Delta}=\xi_{(\tau_b +\Delta)-}-\alpha_1 = c \quad
\text{and} \quad  \xi_{S_1}=\xi_{S_1-}+\beta_1=q.
\]
Here, for any $x \in \mathbb{R}_{+}$, $\tau_x \triangleq \inf\{t
\geq 0: \xi^0_t=x\}$. Let us introduce
\begin{equation}
u_0(\xi) \triangleq \mathbb{E}^{\xi}_0[e^{(b-r)\tau_d}1_{\{\tau_d<
\tau_p\}} u(d)]+\mathbb{E}^{\xi}_0[e^{(b-r)\tau_p} 1_{\{\tau_d>
\tau_p\}} u(p)],
\end{equation}
in which
\begin{equation}\label{eq:ud}
u(d)= \mathbb{E}^{d}_0\left[e^{(b-r) \Delta}
(\bar{C}_{1}(\xi^0_{\Delta-},c)+u(c)) \right] \quad \text{and} \quad
u(p)=\bar{C}_2(p,q)+u(q).
\end{equation}
From (\ref{eq:u})-(\ref{eq:ud}) it can be seen that
\begin{equation}\label{eq:arbitrary}
u(\xi) =
\begin{cases}
\bar{C}_2(\xi,q)+u_0(q), & \xi\leq p; \\ u_0(\xi), & p\leq \xi
\leq d;
\\ r(\xi,c)+ e^{(b-r)\Delta} u_0(c), & \xi \geq d.
\end{cases}
\end{equation}
in which
\begin{equation}\label{eq:r}
r(\xi,c) \triangleq \mathbb{E}^{\xi}_0\left[e^{(b-r) \Delta}
\bar{C}_1 (\xi^0_{\Delta-},c) \right].
\end{equation}
Let us denote the fundamental solutions of
$(\mathcal{A}+(b-r))f=0$, by $\psi$ (increasing) and $\varphi$
(decreasing), and introduce $F \triangleq \psi/\varphi$. Using
(\ref{eq:laplace}), on the interval $(p,d)$ we can write $u$ as
\begin{equation}\label{eq:u-labor}
\frac{u(\xi)}{\varphi(\xi)}=\frac{u(d)}{\varphi(d)}
\frac{(F(\xi)-F(p))}{ (F(d)-F(p))}+\frac{u(p)}{\varphi(p)}\frac{
(F(d)-F(\xi))}{ (F(d)-F(p))}, \quad \xi \in (p,d).
\end{equation}
Then, $W \triangleq \frac{u}{\varphi} \circ F^{-1}$, satisfies
\begin{equation}\label{eq:WF}
W(y)= W(F(d)) \frac{y-F(p)}{F(d)-F(p)}+ W(F(p)) \frac{ (F(d)-y)}{
(F(d)-F(p))}, \quad y \in [F(p),F(d)].
\end{equation}
Using the linear characterization (in the continuation region) of
the band policies in (\ref{eq:WF}),  the following algorithm first
determines the function $u$ for an arbitrary band policy and goes
onto finding the best band policy.

First, let us define
\begin{equation}\label{eq:r1-r2}
R_{1}(x;c) \triangleq \frac{r(\cdot,c)}{\varphi(\cdot)}\circ
F^{-1}(x) \quad \text{and} \quad R_{2}(x;q) \triangleq
\frac{\bar{C}_2(\cdot,q)}{\varphi(\cdot)}\circ F^{-1}(x).
\end{equation}

\noindent {\textbf{\underline{Algorithm}}:}

\begin{enumerate}
\item For a given band policy which is characterized by the
quadruplet $(p,q,c,d)$ such that $p<q<c<d$, we can find the value
function $u$ in (\ref{eq:arbitrary}) using the linear
characterization in (\ref{eq:WF}). On [F(p), F(d)] we will find $W(y)=
\rho y+ \tau$ (in which the slope $\rho$ and the intercept $\tau$
are to be determined) from
\begin{equation}\label{eq:given-q-a}
\begin{split}
e^{(b-r)\Delta}(\rho F(c)+\tau)\frac{\varphi(c)}{\varphi(d)} +R_1(F(d);c)&=\rho F(d)+ \tau,
\\ (\rho F(q)+\tau)\frac{\varphi(q)}{\varphi(p)} +R_{2}(F(p);q)&=\rho
F(p)+\tau.
\end{split}
\end{equation}
$\rho$ and $\tau$ are determined as
\begin{equation}\label{eq:rho-tau}
\begin{split}
\rho=& \frac{\frac{R_2(F(p);q)}{1-\varphi(q)/\varphi(p)}\left(e^{(b-r)\Delta} \frac{\varphi(c)}{\varphi(d)}-1\right)
+R_1(F(d);c)}{F(d)-e^{(b-r) \Delta}\frac{\varphi(c)}{\varphi(d)} F(c)+\frac{\varphi(q)/\varphi(p)\,\, F(q)- F(p)}{1-\varphi(q)/\varphi(p)}\left(
1-e^{(b-r)\Delta} \frac{\varphi(c)}{\varphi(d)}\right)}\, ,
\\ \tau = & \frac{\rho
\left(\frac{\varphi(q)}{\varphi(p)}F(q)-F(p)\right)+R_2(F(p;q))}{1-\frac{\varphi(q)}{\varphi(p)}}\, .
\end{split}
\end{equation}

Now $u$ can be written as
\begin{equation}\label{eq:u-p-q-c-d}
u(\xi)=\begin{cases} u_0(q)+r_{2}(\xi,q), &
x\leq p,\\ u_0(\xi)\triangleq\rho \psi(\xi)+\tau \varphi(\xi), & p \leq x \leq d,
\\ e^{(b-r)\Delta} u_0(c)+r_1(\xi,c), & x
\geq d.
\end{cases}
\end{equation}
From this last expression, we observe that
$(\mathcal{A}+(b-r))u(\xi)=0$ for $\xi \in (p,d)$.

\item
Note that $\rho$ and $\tau$ are functions of $(p,d)$ parametrized
by $(q,c)$. We will find an optimal pair $(p,d)$ given $(q,c)$ by
equating the gradient of the function $(\rho,\tau)$ with respect
to $(p,d)$ to be zero. Now, differentiating the first equation in
(\ref{eq:given-q-a}) with respect to d, and the second with
respect to $p$, and evaluating them at $\tau_d=\rho_d=\tau_p=\rho_p=0$ we obtain
\begin{equation}\label{eq:second-stage}
\begin{split}
&-(\rho F(q)+\tau)\frac{\varphi(q)}{\varphi(p)^2} \varphi'(p)-\rho F'(p)+\frac{\partial}{\partial
y}R_2(y;q)\bigg|_{y=F(p)}F'(p)=0
\\& -e^{(b-r)\Delta}(\rho
F(c)+\tau)\frac{\varphi(c)}{\varphi(d)^2} \varphi'(d) -\rho F'(d)+\frac{\partial}{\partial
y}R_1(y;c) \bigg|_{y=F(d)} F'(d)=0,
\end{split}
\end{equation}
in which $\rho$ and $\tau$ are given by (\ref{eq:rho-tau}). To
find the optimal $(p,d)$ (given $(c,q)$) we solve the
non-linear and implicit system of equations in (\ref{eq:second-stage}).
\begin{remark}\label{rem:smotth-fit-labor}
On $[F(0), F(p)]$ the function $W$ is given by
\begin{equation}\label{eq:R2-shifted}
W(x)=\left((\rho F(q)+\tau)\frac{\varphi(q)}{\varphi(F^{-1}(x))}\right)
+R_{2}(x;q),
\end{equation}
and its left derivative at F(p), $W'(F(p)-)$, is given by
\begin{equation}\label{eq:smooth-fit-at-p}
W'(F(p)-)=-(\rho F(q)+\tau)\frac{\varphi(q)}{\varphi(p)^2} \frac{\varphi'(p)}{F'(p)}+\frac{\partial}{\partial
y}R_2(y;q)\bigg|_{y=F(p)}
\end{equation}
Therefore, the equation in (\ref{eq:second-stage}) in fact implies that the left
and the right derivative of $W$ at $F(p)$ are equal (smooth fit). (Recall that
$W(x)=\rho x+\tau y$ on $[F(p),F(d)]$.) Similarly, the second equation in
(\ref{eq:second-stage}) implies that the left
and the right derivative of $W$ at $F(d)$ are equal. This can be
also expressed as: "$R_2$ shifted by an appropriate amount is
tangential to the line $l(y)=\rho y+\tau$" at $F(p)$.
\end{remark}

\item Next, we vary $q$ and $c$ to find the best
band policy. Such a search can easily carried out in Mathematica.
\end{enumerate}
 To obtain an explicit expression for $g$ in
(\ref{eq:g-labor}) and $r$ in (\ref{eq:r}) we make the following
assumption. We will assume that $\sigma_t=\sigma>0$ (a constant)
in (\ref{eq:demand}).  Now, we can obtain $g$ in
(\ref{eq:g-labor}) (see Appendix) explicitly as
\begin{equation}\label{eq:g-explicit-labor}
g(\xi)= \frac{A^{\mu}}{r-b+(b+\delta) \mu+ \frac{1}{2} \sigma^2
\mu- \frac{1}{2} \sigma^2 \mu^2 } \xi^{\mu}- \frac{w}{r+\delta}
\xi \equiv k_1 \xi^{\mu} + k_2 \xi.
\end{equation}
Note that if $r<b$, then $g(\xi)=\infty$, which implies that
taking no action is optimal. The assumption in Proposition
\ref{prop:existence-uniqueness-labor} that $\max(c_1-c_2,
c_3+c_4)<|k_2|$ is for technical reasons, however it is not very
restrictive. $k_2$ denotes the present value of the total wage
that a firm pays per unit of marginal revenue product of labor and
it should be greater than costs associated with one time hiring or
firing of one unit of marginal  revenue product of labor. Using
(\ref{eq:g-explicit-labor}) we can also calculate $r$ in
(\ref{eq:r}) explicitly as (see Appendix)
\begin{equation}\label{eq:labor-r-exp}
\begin{split}
r(\xi,c) &= e^{(b-r) \Delta} \bigg[-(c_3+c_4)
e^{-(b+\delta)\Delta} \xi
  N(d_1)+ (c_1-c_2) e^{-(b+\delta)\Delta} \xi N(-d_1)
  \\&+ c_3 c N(d_2)-c_1 c N(-d_2) -k_1 \exp\left( \epsilon \right) \xi^{\mu}  -k_2 e^{-(b+\delta)\Delta} \xi+k_1 c^{\mu}+k_2 c \bigg]
\end{split}
\end{equation}
in which
\begin{equation}\label{eq:d1-d2-eps}
\begin{split}
d_1 &\triangleq \frac{1}{\sigma \sqrt{\Delta}}
\log\left(\frac{\xi}{c}\right)+\left(\frac{1}{2}\sigma^2-(b+\delta)\right)\frac{\sqrt{\Delta}}{\sigma},
\\ d_2 &\triangleq \frac{1}{\sigma \sqrt{\Delta}}
\log\left(\frac{\xi}{c}\right)-\left(\frac{1}{2}\sigma^2+(b+\delta)\right)\frac{\sqrt{\Delta}}{\sigma},
\\ \epsilon & \triangleq - \left(b+\delta+\frac{1}{2}\sigma^2
(1-\mu)\right)\mu \Delta.
\end{split}
\end{equation}
Here the function $x \rightarrow N(x)$, $x \in \mathbb{R}$, denotes the cumulative distribution function of an
$N(0,1)$ (standard Gaussian) random variable.
 The infinitesimal generator $\mathcal{A}$ of the process
$\xi$ is $\mathcal{A}u(x) \triangleq (\sigma^2/2) x^2
u''(x)-(b+\delta) x u'(x)$, acting on smooth test functions
$u(\cdot)$. Therefore the fundamental solutions of the equation
$(\mathcal{A}+(b-r)) u=0$ are
\begin{equation}
\psi(x)\triangleq x^{\beta_1}, \quad \varphi \triangleq
x^{\beta_2},
\end{equation}
in which $\beta_1>1$ and $\beta_2<0$ are the roots of the
following quadratic equation (in terms of $\beta$)
\begin{equation}
\frac{1}{2}\sigma^2 \beta^2 - \left(\frac{1}{2} \sigma^2 +
(b+\delta)\right)\beta+ b-r=0.
\end{equation}

The next proposition justifies the second stage of our algorithm.
\begin{proposition}\label{prop:existence-uniqueness-labor}
For a given $(q,c) \in \mathbb{R}^2$, such that $
(c_1q-(k_1q^\mu+k_2q))<0$ there exits a unique solution
$(p^*,d^*)$ to the system of equations (\ref{eq:second-stage}) if
we further assume that $\max(c_1-c_2, c_3+c_4)<|k_2|$. Moreover,
$u^{p^*,q,c,d^*}(x)=\sup_{0<p<d}u^{p,q,c,s}(x), \quad x \geq 0$.
\end{proposition}
\begin{proof}
The proof is similar to that of
Proposition~\ref{prop:forex-example}. Also, see the remark below.
\end{proof}

\begin{remark}\label{rem:extension}
The proof of Proposition~\ref{prop:existence-uniqueness-labor}
only relies on the following properties of the functions $R_1$ and
$R_2$ defined in (\ref{eq:r1-r2}): 1) There exists a point $j \in
(0,\infty)$ such that $y \rightarrow R_1(y;c)$ is concave and
increasing on $(j,\infty)$; 2) $\lim_{y \rightarrow
\infty}R_1(y;c)=\infty$; 3) The function $y \rightarrow R_2(y,q)$
is increasing and concave on $(0,t)$ for some $t<F(q)$ and
decreasing on $(t,\infty)$; 4) Both $y \rightarrow R_1(y;c)$ and
$y \rightarrow R_2(y,q)$ are differentiable.
\end{remark}

Our results in this section can be generalized to the two-sided
control of any one-dimensional diffusion and penalty functions
satisfying the conditions in Remark~\ref{rem:extension} are
satisfied. It is worth pointing out that Weeransinghe \cite{wee}
has studied the two-sided bounded variation control within the
framework of \emph{singular stochastic control} of linear
diffusions for a large class of cost functions by using of the
functional relationship between the value function of optimal
stopping and that of singular stochastic control (see e.g.
Karatzas and Shreve \cite{Karat-Shr1985}).

\begin{figure}\label{figure:2}
\begin{center}
\begin{minipage}{0.45\textwidth}
\centering \includegraphics[scale=0.75]{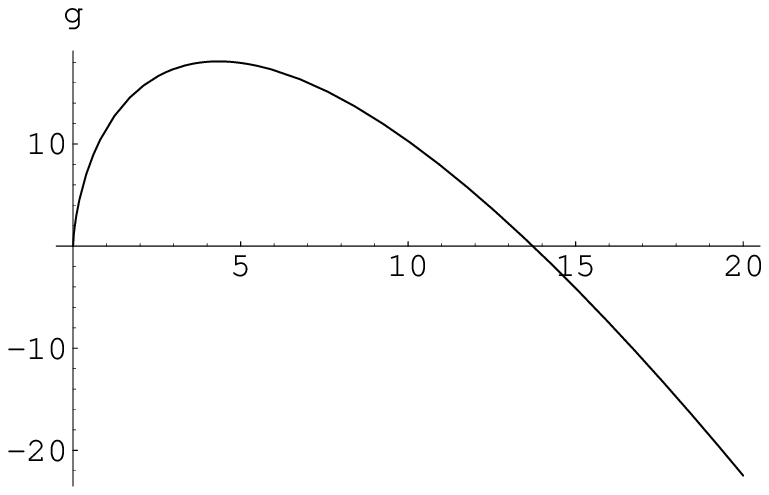} \\
(a)
\end{minipage}
\begin{minipage}{0.45\textwidth}
\centering \includegraphics[scale=0.75]{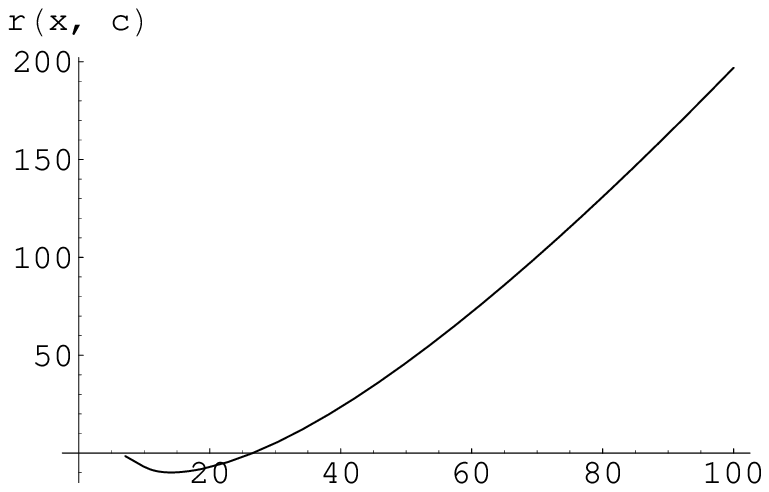} \\
(b)
\end{minipage}
\begin{minipage}{0.45\textwidth}
\centering \includegraphics[scale=0.75]{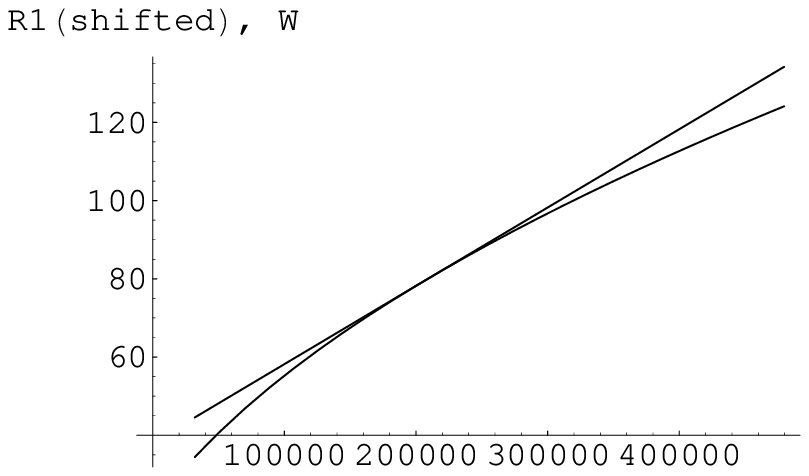} \\
(c)
\end{minipage}
\begin{minipage}{0.45\textwidth}
\centering \includegraphics[scale=0.75]{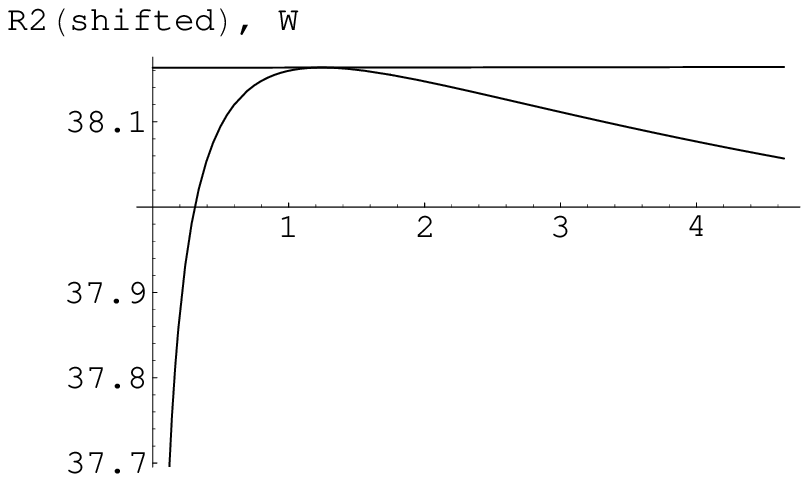} \\
(d)
\end{minipage}
\begin{minipage}{0.45\textwidth}
\centering \includegraphics[scale=0.75]{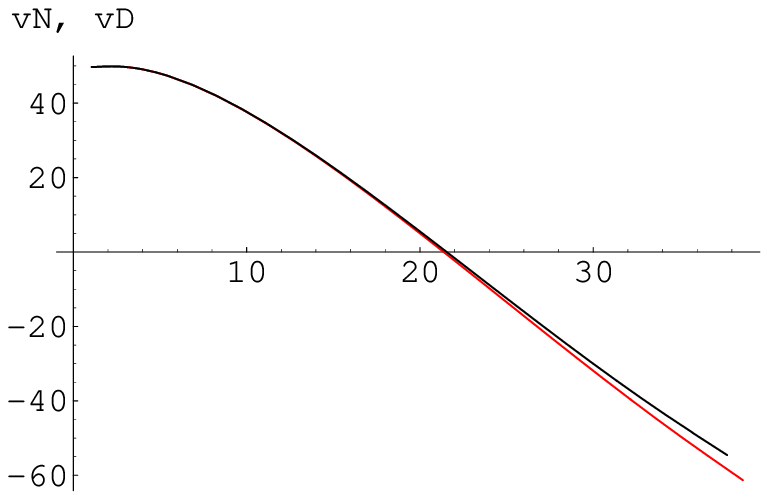} \\
(e)
\end{minipage}
\begin{minipage}{0.45\textwidth}
\centering \includegraphics[scale=0.75]{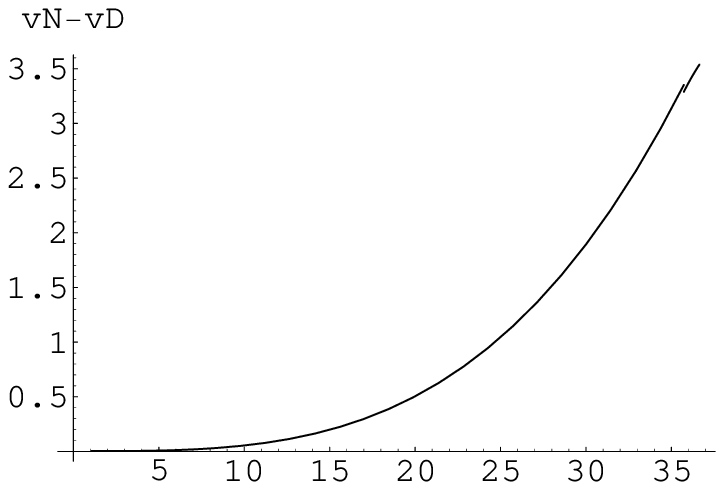} \\
(f)
\end{minipage}
\begin{minipage}{0.45\textwidth}
\centering \includegraphics[scale=0.75]{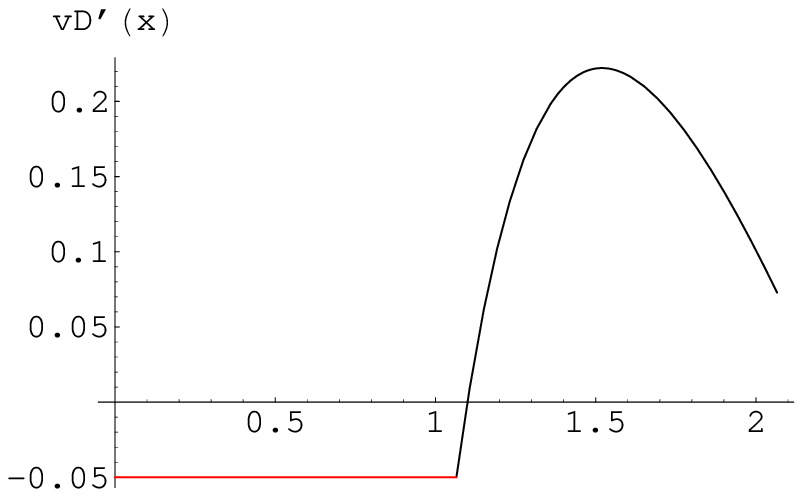} \\
(g)
\end{minipage}
\begin{minipage}{0.45\textwidth}
\centering \includegraphics[scale=0.75]{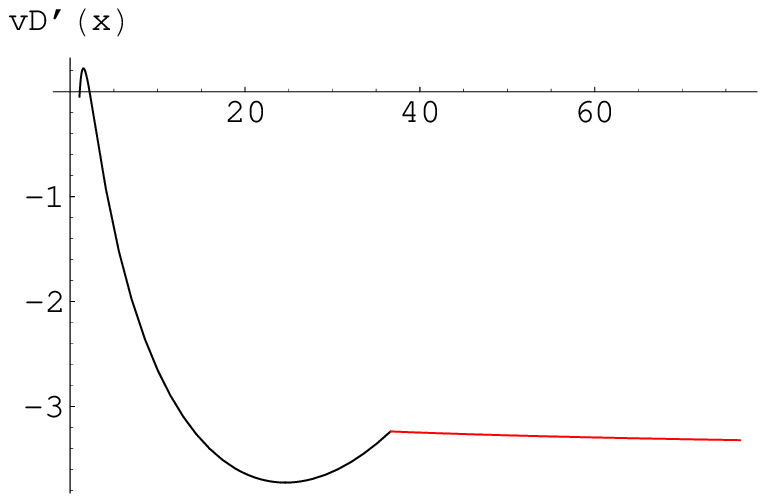} \\
(h)
\end{minipage}
\caption{\small (a) The graph of $g(x)$. (b) The graph of $r(x,
c^*)$ for $\Delta>0$ (c) The graph of line $\rho^*y+\tau^*$ we
obtain via our algorithm and $R_1(y, c^*)$ after it is shifted
vertically by $e^{(b-r)\Delta}(\rho
F(c)+\tau)\frac{\varphi(c)}{\varphi(F^{-1}(y))}$. (d) The graph of
the line $\rho^* y+\tau^*$ and $R_2(y, c^*)$ after it is shifted
(see (\ref{eq:R2-shifted}) for the amount of shift).
 (e) The two value functions,  $v^N(x)$ ($\Delta=0$) above and $v^D(x)$ ($\Delta>0$) below.
  (f) Plot of difference, $v^N(x)-v^D(x)$.
   (f) Plot of difference, $v^N(x)-v^D(x)$. (g) (h)The derivatives match at $x=p$ and $x=d$ ($\Delta>0)$.}
\end{center}

\end{figure}

\subsection{Numerical Example}
In this section, we will give a numerical example for the labor problem with and without delay.
We select
the parameters as $b=0.03, r=0.06, \mu=0.75, \sigma=0.35, \delta=0.1,
A=5, w=2, \Delta=0.5$, $c_1=0.05$,
$c_2=0.1$, $c_3=2$ and $c_4=1$.  The results
we obtain are summarized in the following table:
\begin{center}
\begin{tabular}{|c|c|c|c|c|c|c|}
\hline & $\rho$ & $\tau$ & $p$ & $q$ & $c$ & $d$ \\ \hline
$\Delta=0$ & $0.0002003$ & $38.1633$ & $1.0664$ & 2.125 & 7.240 &
35.728 \\ \hline $\Delta=0.5$ & 0.0001725 & 38.1597 &1.0661 &
2.100 & 7.120 & 36.640 \\ \hline
\end{tabular}
\end{center}
Both the slope $\rho$ and the intercept $\tau$ are greater in the
no-delay case and therefore, the value function corresponding to
no-delay problem $v^N(x)$ will dominate that to delay problem
$v^D(x)$.  On the right boundary, we have $(7.240, 35.728)\subset
(7.120, 36.640)$ and on the left boundary $(c, d)$ pair has
shifted to the left with delay. As a result, the continuation
region $(p, d)$ has expanded with delay: $\mathbf{C}^N\triangleq
(1.0664, 35.728)\subset (1.0661, 36.640)\triangleq \mathbf{C}^D$.
An explanations for this phenomenon can be made through the
relative size of costs of firing and hiring, the size of delay
parameter, the shape of $g$ function, etc.  In our example, the
firing cost is relatively larger than hiring cost, the penalty of
firing becomes smaller with delay (than without delay) which
encourages the controller not make hasty firing decisions, facing
relatively large firing costs. Or since there is a chance that the
process moves to the left during the delay period due to voluntary
quits, this effect may help to reduce firing costs even though the
decision making is postponed.

\section{Conclusion}\label{sec:conclusion}
In this paper we give a new characterization of the \emph{value
function} of one-sided and two-sided impulse control problems with
implementation delays. We also provided easily implemented
algorithms to find out the optimal control and the value function.
Our methodology bypasses the need to guess the form of solution of
quasi-variational inequalities and prove that this solution
satisfies a verification lemma.   Since our method directly finds
the value function, we believe that this method can solve a larger
set of problems than just with quasi-variational inequalities.
Indeed, we applied our results to solving some specific examples.
As an important application of a two-sided impulse control problem
with decision delays we found out the optimal hiring and firing
decisions of a firm facing regulatory delays and stochastic
demand.

Here we considered a problem in which the decision maker needs to decide
whether to take action and, after some delay, needs to decide the
magnitude of her action. In the future, we will consider problems in which the decision maker takes action and
waits that action to be implemented. We will also consider a general characterization of the value
function and the optimal controls when the decision delay is not a
constant but it depends on the magnitude of the action taken as in
\cite{sub-jarrow} or it depends on the value of the state variable
that is controlled as in \cite{alvarez-keppo}.

\noindent\textbf{Acknowledgment}

We are grateful to the the referee for his/her detailed comments
that helped us improve the manuscript.
\section{Appendix}
\subsection{Derivations of (\ref{eq:g-labor}) and (\ref{eq:labor-r-exp})}
Using (\ref{eq:xi0t}) we can write (\ref{eq:g-labor}) as
\begin{equation}
\begin{split}
g(\xi)&=\mathbb{E}^{\xi}_0\left[\int_0^{\infty}A^{\mu} \xi^{\mu} e^{(b-r)t} \exp (-(b+\delta)\mu
t- \sigma \mu B_t -\frac{1}{2} \sigma^2 \mu t) dt\right]
\\ & \qquad \qquad \qquad \qquad \qquad \qquad - w \mathbb{E}_{0}^{\xi}\left[\int_0^{\infty} \xi \exp (-(b+\delta)
t- \sigma B_t -\frac{1}{2} \sigma^2 t)dt \right]
\\&=A^{\mu} \xi^{\mu}
\int_0^{\infty}\exp\left[t\left(b-r-(b+\delta)\mu-\frac{1}{2}\sigma^2 \mu +\frac{1}{2}\sigma^2 \mu^2\right)\right]dt
-w \xi \int_0^{\infty} \exp(-(b+\delta)t) dt ,
\end{split}
\end{equation}
from which we obtain (\ref{eq:g-explicit-labor}) under the assumption that $r>b$. Here the second inequality follows from the Fubini's theorem and using the Laplace transform of $B_t$.

In what follows we will present the derivation of
(\ref{eq:labor-r-exp}). We can write (\ref{eq:r}) as
    \begin{equation}\label{eq:r-appendix}
\begin{split}
r(\xi,c)&=e^{(b-r)\Delta}\mathbb{E}_0^{\xi}\bigg[\left(-c_3(\xi_{\Delta}-c)-c_4 \xi_{\Delta}\right)1_{\{\xi_{\Delta}>c\}}
+(-c_1(c-\xi_{\Delta})-c_2 \xi_{\Delta}) 1_{\{\xi_{\Delta}<c\}}
\\& -k_1 \xi_{\Delta}^{\mu}-k_2 \xi_{\Delta}+k_1 c^{\mu}+k_2 c\bigg]
\end{split}
\end{equation}

Using (\ref{eq:xi0t}) and the assumption that $\sigma_t=\sigma \in
\mathbb{R}_{+}$, we compute
\begin{equation}\label{eq:A-B-C}
\begin{split}
A &\triangleq
\mathbb{E}^{\xi}_{0}\left[1_{\{\xi_{\Delta}>c\}}\right]=N(d_2),
\quad B \triangleq
\mathbb{E}^{\xi}_{0}\left[1_{\{\xi_{\Delta}<c\}}\right]=1-A=N(-d_2),
\\ C(\theta) &\triangleq \mathbb{E}^{\xi}\left[\xi_{\Delta}^{\theta}\right]
= \xi^{\theta}\exp\left(-\left(b+\delta+\frac{1}{2}\sigma^2(1-\theta)\right)\theta
\Delta\right),
\end{split}
\end{equation}
where $\theta=1$ or $\theta=\mu$. Here the third equality follows from the Laplace transform of $B_t$
We will also need to compute
\begin{equation}\label{eq:defn-of-D}
D \triangleq
\mathbb{E}^{\xi}_0\left[\xi_{\Delta}1_{\{\xi_{\Delta}>c\}}\right].
\end{equation}
We will denote
\[
\kappa=\exp\left(-\frac{1}{2} \sigma^2 \Delta+ \sigma \sqrt{\Delta} \eta
\right),
\]
in which $\eta=B_{\Delta}/\sqrt{\Delta}$, is an $N(0,1)$ random variable. Then $\xi_{\Delta}= \xi \exp(-(b+\delta)\Delta)
\kappa$ and $A=\xi e^{-(b+\delta)\Delta} \mathbb{E}^{\xi}_{0}\left[1_{\{\xi_{\Delta}>c\}}
\kappa\right]$. Introducing a new probability measure Q by the radon-nikodym derivative
$dQ^{\xi}/dP^{\xi}_{0}=\kappa$, we get
\[
D=e^{-(b+\delta)\Delta} \xi Q^{\xi}(\xi_{\Delta}>c).
\]
Under the measure $Q^{\xi}$, $n \triangleq -\eta -\sigma \sqrt{\Delta}$ is
$N(0,1)$ and we can write $\xi_{\Delta}$ in terms of $n$ as
\begin{equation}\label{eq:xi-under-Q}
\xi_{\Delta}=\xi \exp\left(-(b+\delta-\frac{1}{2}\sigma^2)\Delta+ \sigma
\sqrt{\Delta}n\right).
\end{equation}
Using (\ref{eq:xi-under-Q}), we can compute
\begin{equation}\label{eq:sol-A}
D= \xi e^{-(b+\delta)\Delta}N(d_1),
\end{equation}
in which $d_1$ is given by (\ref{eq:d1-d2-eps}). We can then
immediately obtain,
\begin{equation}\label{eq:defn-E}
E \triangleq
\mathbb{E}^{\xi}_0\left[\xi_{\Delta}1_{\{\xi_{\Delta}<c\}}\right]=\xi e^{-(b+\delta)\Delta}(1-Q^{\xi}(\xi_{\Delta}>c))=\xi
e^{-(b+\delta)\Delta}N(-d_1).
\end{equation}
Using (\ref{eq:r-appendix}), (\ref{eq:A-B-C}),
(\ref{eq:defn-of-D}) and (\ref{eq:defn-E}) we obtain (\ref{eq:r}).

\subsection{A Technical Lemma}

\begin{lemma}\label{lem:Lipschitz}
Define
\begin{equation*}
G(x,\gamma) \triangleq \sup_{\tau \in S}\mathbb{E}^{x}[e^{-\alpha
\tau}(h(X^0_\tau)+\gamma e^{-\alpha \Delta})], \quad x \in \mathbb{R}, \, \gamma \in
\mathbb{R},
\end{equation*}
for some Borel function $h$. Then for $\gamma_1>\gamma_2$ we have
that
\[
G(x,\gamma_1)-G(x,\gamma_2) \leq \gamma_1-\gamma_2.
\]
\end{lemma}
\begin{proof}
See the proof of Lemma 3.3 in \cite{DE2005}.
\end{proof}

\subsection{Proof of Proposition ~\ref{prop:forex-example}}
The proof follows from the analysis of the function $r$. The
following remark will be helpful in the analysis that follows.

\begin{remark}\normalfont
Let us denote $H(y)\triangleq(h/\varphi)\circ(F^{-1}(y)), y>0$. If
$h(\cdot)$ is twice-differentiable at $x\in \mathcal{I}$ and
$y\triangleq F(x)$, then $H^{'}(y)=m(x)$ and
$H^{''}(y)=m^{'}(x)/F^{'}(x)$ with
\begin{equation}\label{eq:devH}
m(x)= \frac{1}{F^{'}(x)}\left(\frac{h}{\varphi}\right)^{'}(x),
\quad \text{and} \quad H^{''}(y) [(\mathcal{A}-\alpha)h(x)]\geq 0,
\quad y=F(x),
\end{equation}
with strict inequality if $H^{''}(y)\neq 0$.
\end{remark}
\subsubsection{The Analysis of the Function $r$ in (\ref{eq:r-forex-example})}

 Let us check the sign of
$\left(\frac{r}{\varphi}\right)'(x)=\frac{r'\varphi-r\varphi'}{\varphi^2}(x)$
which is the same as the derivative of $R$ as can be observed from
the first equation in (\ref{eq:devH}).   The sign of
$\left(\frac{r}{\varphi}\right)'(x)$ is the same as that of
\begin{align}\label{eq:check}
&\frac{\sqrt{2\alpha}}{\alpha}\left(x^2-a^2+\Delta-2\alpha\lambda\Delta \exp\left(-\frac{(a-x)^2}{4\Delta^2}\right)-c\alpha\right)\nonumber\\
&+\lambda(a-x)\left(-\frac{1}{\Delta} \exp\left(-\frac{(a-x)^2}{4\Delta^2}\right)+\frac{1}{\Delta}\phi\left(\frac{a-x}{\Delta}\right)+\sqrt{2\alpha}\left(2N\left(\frac{a-x}{\Delta}\right)-1\right)\right)\nonumber\\
&+\frac{2x}{\alpha}+\lambda\left(2N\left(\frac{a-x}{\Delta}\right)-1\right).
\end{align}
Using the fact $2N\left(\frac{a-x}{\Delta}\right)<1$ for $x>a$ and
$-\frac{1}{\Delta}
\exp\left(-\frac{(a-x)^2}{4\Delta^2}\right)+\frac{1}{\Delta}\phi\left(\frac{a-x}{\Delta}\right)<0$
for $x>a$ sufficiently large,in this equation (for sufficiently
large x)  we identify the absolute value of the negative terms as
$\frac{\sqrt{2\alpha}}{\alpha}\lambda\Delta
\exp\left(-\frac{(a-x)^2}{4\Delta^2}\right)<\frac{\sqrt{2\alpha}}{\alpha}\lambda\Delta$
, $c\alpha$ and
$|\lambda\left(2N\left(\frac{a-x}{\Delta}\right)-1\right)|<\lambda$.
Since these negative terms are bounded, if we take sufficiently
large value, say $a'$, the sign of (\ref{eq:check}) is positive
for $x\in (a', \infty)$.  Moreover, we can directly calculate
$\lim_{y\rightarrow +\infty}\frac{\partial}{\partial y}R(y;a)=0$
to check the behavior of $R(y;a)$ for a large $y$. We also know
that $R(y;a)\triangleq (r(\cdot, a)/\varphi(\cdot))\circ
F^{-1}(y)$ is negative at $y=F(a)$. On the other hand,
$1/\varphi(F^{-1}(y))=\sqrt{y}$ is increasing and concave
function. It follows that $R(y;
c)+\frac{\gamma}{\varphi(F^{-1}(y))}$ is an increasing function on
$y\in (F(a'), \infty)$.

To investigate the concavity of $R(y; a)$, we set
\begin{align*}
q(x, a)&\triangleq
\frac{1}{2}x^2\frac{\lambda}{\Delta}\left(e^{-\frac{(a-x)^2}{\Delta}}\left(1-\frac{2(a-x)^2}{4\Delta^2}\right)
-3\phi\left(\frac{a-x}{\Delta}\right)-\lambda(a-x)\phi'\left(\frac{a-x}{\Delta}\right)\right)\\
&+\alpha x^2 -\alpha r(x, a)
\end{align*}
so that $(\mathcal{A}-\alpha)r(x, a)=q(x, a)$ for every $x>0$. We
have $\lim_{x\rightarrow \infty}q(x)=-\infty$ if $\alpha<4$.  By
the second equation in (\ref{eq:devH}), the function $R(y;a)$
becomes concave eventually. Since $R(\cdot;a)$ is increasing and
concave on $(a'',\infty)$ for some $a''>a'$ and $\lim_{y
\rightarrow \infty}R(y;a)=\infty$  we can find a unique linear
majorant to $R^{\gamma}(\cdot,a)$ in
Lemma~\ref{lem:majorant-forex} (the linear majorant majorizes
$R^{\gamma}(\cdot,a)$ in the continuation region and is equal to
$R^{\gamma}(\cdot,a)$ in the stopping region). The rest of the
proof from Proposition ~\ref{prop:existence-uniquenes-forex}.

\bibliographystyle{plain}

{\small

\begin{thebibliography}{10}

\bibitem{alvarez2}
L.~H.~R. Alvarez.
\newblock A class of solvable impulse control problems.
\newblock {\em Appl. Math. and Optim.}, 49:265--295, 2004.

\bibitem{alvarez1}
L.~H.~R. Alvarez.
\newblock Stochastic forest stand value and optimal timber harvesting.
\newblock {\em SIAM J. Control. Optim.}, 42(6):1972--1993, 2004.

\bibitem{alvarez-keppo}
L.~H.~R. Alvarez and J.~Keppo.
\newblock The impact of delivery lags on irreversible investment under
  uncertainty.
\newblock {\em European Journal of Operational Research}, 136:173--180, 2002.

\bibitem{alvarez3}
L.~H.~R. Alvarez and J.~Virtanen.
\newblock A class of solvable stochastic dividend optimization problems: On the
  general impact of flexibility on valuation.
\newblock {\em Economic Theory}, 28:373--398, 2006.

\bibitem{bar-ilan-perry-stadje}
A.~Bar-Ilan, D.~Perry, and W.~Stadje.
\newblock A generalized impulse control of cash management.
\newblock {\em Journal of Economic Dynamics and Control}, 28:1013--1033, 2004.

\bibitem{bar-ilan-strange1}
A.~Bar-Ilan and W.~C. Strange.
\newblock Investment lags.
\newblock {\em American Economic Review}, 86:610--622, 1996.

\bibitem{bar-ilan}
A.~Bar-Ilan and W.~C. Strange.
\newblock A model of sequential invetment.
\newblock {\em Journal of Economic Dynamics and Control}, 22:437--463, 1998.

\bibitem{ben-lions}
A.~Bensoussan and J.~L. Lions.
\newblock {\em Impulse Control and Quasi-Variational Inequalities}.
\newblock Gauthier-Villars, Paris, 1982.

\bibitem{bb}
S.~Bentolila and G.~Bertola.
\newblock Firing costs and labor demand: How bad is {E}urosclerosis?
\newblock {\em Review of Economic Studies}, 57:381--402, 1990.

\bibitem{salminen}
A.~N. Borodin and P.~Salminen.
\newblock {\em Handbook of Brownian Motion Facts and Formulae}.
\newblock Birkh\"{a}user, Boston, 2002.

\bibitem{DE2005}
S.~Dayanik and M.~Egami.
\newblock Solving stochastic impulse control problems via optimal stopping for
  one-dimensional diffusions.
\newblock {\em preprint, www.umich.edu/$\sim$egami}, 2005.

\bibitem{DK2003}
S.~Dayanik and I.~Karatzas.
\newblock On the optimal stopping problem for one-dimensional diffusions.
\newblock {\em Stochastic Processes and their Applications}, 107 (2):173--212,
  2003.

\bibitem{dynkin2}
E.~Dynkin.
\newblock Optimal choice of stopping moment of a {M}arkov process.
\newblock {\em Dokl. Akad. Nauk. SSSR}, 150:238--240, 1963.

\bibitem{dynkin}
E.~Dynkin.
\newblock {\em Markov processes, Volume II}.
\newblock Springer Verlag, Berlin, 1965.

\bibitem{elsanosi-oks-sulem}
I.~Elsanosi, B.~{\O}ksendal, and A.~Sulem.
\newblock Some solvable stochastic control problems with delay.
\newblock {\em Stochastics and Stochastics Reports}, 71:225--243, 2000.

\bibitem{Karat-Shr1985}
I.~Karatzas and S.~E. Shreve.
\newblock Connections between optimal stopping and singular stochastic control
  ii. reflected follower problems.
\newblock {\em SIAM J. Control Optim.}, 23 (3):433--451, 1985.

\bibitem{keppo-peura}
J.~Keppo and S.~Peura.
\newblock Optimal bank capital with costly recapitalization.
\newblock {\em \emph{To appear} in the Journal of Business}, 2005.

\bibitem{mundaca-oksendal}
G.~Mundaca and B.~{\O}ksendal.
\newblock Optimal stochastic intervention control with application to the
  exchange rate.
\newblock {\em Journal of Mathematical Economics}, 29:225--243, 1998.

\bibitem{O1999}
B.~{\O}ksendal.
\newblock Stochastic control problems where small intervention costs have big
  effects.
\newblock {\em Appl. Math. Optim.}, 40:355--375, 1999.

\bibitem{oksendal0sulem2001}
B.~{\O}ksendal and A.~Sulem.
\newblock {\em A Maximum Principle for optimal control of stochastic systems
  with delay, with applications to finance}, pages 64--79.
\newblock IOS Press, 2001.
\newblock Optima Control and Partial Differential Equations, J.L. Menaldi et.
  al.(editors)).

\bibitem{oksendal-book-2}
B.~{\O}ksendal and A.~Sulem.
\newblock {\em Applied stochastic controll of jump diffusions}.
\newblock Springer-Verlag, New York, 2005.

\bibitem{oksendal-delay-impulse}
B.~{\O}ksendal and A.~Sulem.
\newblock Optimal stochastic impulse control with delayed reaction.
\newblock {\em Preprint. University of Oslo}, 2005.

\bibitem{shepp-shiryaev}
L.~A. Shepp and A.~N. Shiryaev.
\newblock Hiring and firing optimally in a large corporation.
\newblock {\em Journal of Economic Dynamics and Control}, 20:1523--1540, 1996.

\bibitem{sub-jarrow}
A.~Subramanian and R.~A. Jarrow.
\newblock The liquidity discount.
\newblock {\em Mathematical Finance}, 11:447--474, 2001.

\bibitem{wee}
A.~Weerasinghe.
\newblock A bounded variation control problem for diffusion process.
\newblock {\em SIAM J. Control Optim.}, 44(2):389--417, 2005.

\end{thebibliography}

}
\end{document}